\documentclass[english]{article}
\usepackage[T1]{fontenc}
\usepackage[latin9]{inputenc}
\usepackage[letterpaper]{geometry}
\usepackage{float}
\usepackage{amsmath}
\usepackage{graphicx}
\usepackage{color}
\usepackage{epstopdf}
\usepackage{multirow}
\usepackage{array}
\usepackage{mathrsfs}

\newcommand{\average}[1]{\{\!\!\!\{{#1}\}\!\!\!\}}
\newcommand{\jump}[1]{[\![{#1}]\!]}

\newtheorem{theorem}{Theorem}

\newtheorem{lemma}{Lemma}

\usepackage{amsfonts} 
\usepackage{multirow}
\newcommand{\specialcell}[2][c]{\begin{tabular}[#1]{@{}c@{}}#2\end{tabular}}
\makeatother

\usepackage{babel}

\title{Generalized Multiscale Finite Element Method for Elasticity Equations}

\author{Eric T. Chung\thanks{Department of Mathematics, The Chinese University of Hong Kong, Hong Kong SAR. 
This research is partially supported by the Hong Kong RGC General Research Fund (Project number: 400411).},
Yalchin Efendiev\thanks{Department of Mathematics, Texas A\&M University, College Station, TX 77843 and  Numerical Porous Media SRI Center, KAUST, Thuwal, Saudia Arabia}, 
and Shubin Fu\thanks{Department of Mathematics, Texas A\&M University, College Station, TX 77843.}
}

\begin{document}
\maketitle

\begin{abstract}

In this paper, we discuss the application of Generalized Multiscale Finite Element 
Method (GMsFEM) to
elasticity equation in heterogeneous media. Our applications are motivated by elastic wave 
propagation in subsurface where the subsurface properties can be highly heterogeneous
and have high contrast. We present the construction of main ingredients for GMsFEM such as
the snapshot space and offline spaces. The latter is constructed using 
local spectral decomposition
in the snapshot space. The spectral decomposition is based on the analysis which is provided
in the paper. We consider both continuous Galerkin and discontinuous Galerkin coupling of basis
functions. Both approaches have their cons and pros. Continuous Galerkin methods allow 
avoiding penalty parameters though they involve partition of unity functions which can 
alter the properties of multiscale basis functions. On the other hand, discontinuous Galerkin
techniques allow gluing multiscale basis functions without any modifications. Because basis
functions are constructed independently from each other, this approach provides an advantage.
We discuss the use of oversampling techniques that use snapshots in larger regions to
construct the offline space.
We provide numerical results to show that one can accurately approximate the solution
using reduced number of degrees of freedom.

\end{abstract}

\section{Introduction}

Many materials in nature are highly heterogeneous and their properties can vary at different scales.
Direct numerical simulations in such multiscale media are prohibitively expensive and some type
of model reduction is needed.  
Multiscale approaches such as homogenization and numerical homogenization \cite{cao2005iterated,abdulle2006analysis,schroder2014numerical,buck2013multiscale,francfort1986homogenization,oleinik2009mathematical,vinh2011homogenized,liu2009multiscale} have been routinely
used to model macroscopic properties and macroscopic behavior of elastic materials.
These approaches compute the effective material properties based on representative volume simulations.
These properties are further used to solve macroscale equations. In this paper, our goal
is to design multiscale method for elasticity equations in the media when the media properties
do not have scale separation and classical homogenization and numerical homogenization
techniques do not work. 
We are motivated by seismic wave applications when elastic wave propagation
in heterogeneous subsurface formation is studied where the subsurface 
properties can contain vugs,
fractures, and cavities of different sizes. In this paper, we develop multiscale methods for
static problems and present their analysis. 
%We refer to \cite{xxx} (SEG Expanded Abstracts), 
%where we have applied these methods to seismic wave propagation.

In this paper, we design a multiscale model reduction techniques using GMsFEM for 
steady state elasticity equation in heterogeneous media
\begin{equation}
\label{eq:elastic1}
{\partial \over \partial x_i} (c_{ijkl}(x) e_{kl}(u))=f_j(x),
\end{equation}
where $e_{kl}(u)={1\over 2}({\partial u_k \over \partial x_l}+{\partial u_l \over \partial x_k})$
and $c_{ijkl}(x)$ is a multiscale field with a high contrast.
GMsFEM has been studied for a various applications related to flow problems
(see \cite{egh12, eglp13, cel14, Efendiev_LS_MSDG_2013, elms2014}).
 In GMsFEM, we solve equation (\ref{eq:elastic1}) on a coarse grid
where each coarse grid consists of a union of fine-grid blocks. In particular, we 
design
(1) a snapshot space (2) an offline space for each coarse patch. The offline space
consists of multiscale basis functions that are coupled in a global formulation.
In this paper, we consider several choices for snapshot spaces, offline spaces, and global
coupling. The main idea of the snapshot space in each coarse patch is to provide
an exhaustive space where an appropriate spectral decomposition is performed.
This space contains local 
functions that can mimic the global solution behavior in the coarse patch
for all right hand sides or boundary conditions. We consider two choices for
the snapshot space. The first one consists of all fine-grid functions in each coarse patch
and the second one consists of harmonic extensions. Next, we propose a local
spectral decomposition in the snapshot space which allows selecting multiscale basis
functions. This local spectral decomposition is based on the analysis and depends
on the global coupling mechanisms. We consider several choices for the local 
spectral decomposition including oversampling approach where larger domains
are used in the eigenvalue problem. The oversampling technique uses
larger domains to compute snapshot vectors that are more consistent 
with local solution space and thus can have much lower dimension.

To couple multiscale basis functions constructed in the offline space, we consider
two methods, conforming Galerkin (CG) approach and discontinuous Galerkin (DG) approach based
on symmetric interior penalty method for (\ref{eq:elastic1}). 
These approaches are studied for linear elliptic equations in \cite{egh12,Efendiev_JCP_2013}.
Both approaches provide a global coupling
for multiscale basis functions where the solution is sought in the space spanned 
by these multiscale
basis functions. This representation
allows approximating the solution
with a reduced number of degrees of freedom. 
The constructions of the basis functions are different
for continuous Galerkin and discontinuous Galerkin methods as the local
spectral decomposition relies on the analysis. In particular, for continuous Galerkin
approach, we use partition of unity functions and discuss several choices for partition
of unity functions. We provide an analysis of both approaches. The offline space construction
is based on the analysis.

We present numerical results where we study the convergence of continuous and 
discontinuous Galerkin
methods using various snapshot spaces as well as with and without the use of oversampling. 
We consider highly heterogeneous coefficients that contain high contrast.
Our numerical
results show that the proposed approaches allow approximating the solution accurately with a fewer degrees
of freedom. In particular, when using the snapshot space consisting of harmonic extension functions,
we obtain better convergence results. In addition, oversampling methods and the use of snapshot spaces
constructed in the oversampled domains can substantially improve the convergence.

The paper is organized as follows. In Section \ref{sec:prelim}, we state the problem and the notations for coarse and fine grids. In Section \ref{sec:gmsfem}, we give the construction of multiscale basis functions, snapshot spaces
and offline spaces, as well as global coupling via CG and DG. 
In Section \ref{sec:gmsfem_num_res}, we present numerical results.
Sections \ref{sec:gmsfem_error}-\ref{sec:gmsfem_error_DG} are
devoted to the analysis of the methods.

\section{Preliminaries}
\label{sec:prelim}

In this section, we will present the general framework of GMsFEM 
for linear elasticity in high-contrast media. 
Let $D\subset\mathbb{R}^{2}$ (or $R^3$) be a bounded domain representing the elastic body of interest,
and let $ {u} = (u_1,u_2)$ be the displacement field.
The strain tensor $ {\epsilon}( {u}) = (\epsilon_{ij}( {u}))_{1\leq i,j \leq 2}$
is defined by
\begin{equation*}
 {\epsilon}( {u}) = \frac{1}{2} ( \nabla  {u} + \nabla  {u}^T ),
\end{equation*}
where $\displaystyle \nabla  {u} = (\frac{\partial u_i}{\partial x_j})_{1\leq i,j \leq 2}$.
In the component form, we have
\begin{equation*}
\epsilon_{ij}( {u}) = \frac{1}{2} \Big( \frac{\partial u_i}{\partial x_j} + \frac{\partial u_j}{\partial x_i} \Big), \quad 1\leq i,j \leq 2. 
\end{equation*}
In this paper, we assume the medium is isotropic.
Thus, the stress tensor $ {\sigma}( {u}) = (\sigma_{ij}( {u}))_{1\leq i,j \leq 2}$
is related to the strain tensor $ {\epsilon}( {u})$ in the following way
\begin{equation*}
 {\sigma} = 2\mu  {\epsilon} + \lambda \nabla\cdot  {u} \,  {I},
\end{equation*}
where $\lambda>0$ and $\mu>0$ are the Lam\'e coefficients.
We assume that $\lambda$ and $\mu$ have highly heterogeneous spatial variations with high contrasts. 
Given a forcing term $ {f} = (f_1,f_2)$, the displacement field $ {u}$ satisfies
the following
\begin{equation}
\label{ob_equ}
- \nabla \cdot  {\sigma} =  {f}, \quad\text{ in } \; D
\end{equation}
or in component form
\begin{equation}
- \Big(\frac{\partial \sigma_{i1}}{\partial x_1} + \frac{\partial \sigma_{i2}}{\partial x_2} \Big) = f_i, \quad \text{ in } \; D, \quad i=1,2.
\end{equation}
For simplicity, we will consider the homogeneous Dirichlet boundary condition $ {u} =  {0}$ on $\partial D$. 
%Other types of boundary conditions can be taken care easily in the way used in classical approaches. 

Let $\mathcal{T}^H$ be a standard triangulation of the domain $D$
where $H>0$ is the mesh size. We call $\mathcal{T}^H$ the coarse grid 
and $H$ the coarse mesh size. 
Elements of $\mathcal{T}^H$ are called coarse grid blocks. 
The set of all coarse grid edges is denoted by $\mathcal{E}^H$
and the set of all coarse grid nodes is denoted by $\mathcal{S}^H$.
We also use $N_S$ to denote the number of coarse grid nodes, $N$ to denote the number of coarse grid blocks. 
In addition, we let $\mathcal{T}^h$
be a conforming refinement of the triangulation $\mathcal{T}^H$.
We call $\mathcal{T}^h$ the fine grid and $h>0$ is the fine mesh size. 
We remark that the use of the conforming refinement is only to simplify the discussion
of the methodology and is not a restriction of the method.

Let $V^h$ be a finite element space defined on the fine grid. 
The fine-grid solution $ {u}_h$ can be obtained as
\begin{equation}
\label{cg_fine_sol}
a( {u}_h,  {v}) = ( {f},  {v}), \quad \forall  {v}\in V^h,
\end{equation}
where
\begin{equation}
a( {u},  {v}) = \int_D \Big( 2\mu  {\epsilon}( {u}) :  {\epsilon}( {v}) 
+ \lambda \nabla\cdot  {u} \, \nabla\cdot  {v}  \Big) \; d {x},
\quad
( {f},  {v}) = \int_D  {f} \cdot  {v} \; d {x}
\end{equation}
and
\begin{equation}
 {\epsilon}( {u}) :  {\epsilon}( {v}) = \sum_{i,j=1}^2 \epsilon_{ij}( {u}) \epsilon_{ij}( {v}),
\quad
 {f} \cdot  {v} = \sum_{i=1}^2 f_i v_i.
\end{equation}

Now, we present GMsFEM. 
The discussion consists of two main steps, namely,
the construction of local basis functions
and the global coupling. 
In this paper, we will develop and analyze two types of global coupling,
namely, the continuous Galerkin coupling and the discontinuous Galerkin coupling.
These two couplings will require two types of local basis functions. 
In essence, the CG coupling will need vertex-based local basis functions
and the DG coupling will need element-based local basis functions. 

\begin{figure}[ht]
\centering
\includegraphics[width=17cm,height=8cm]{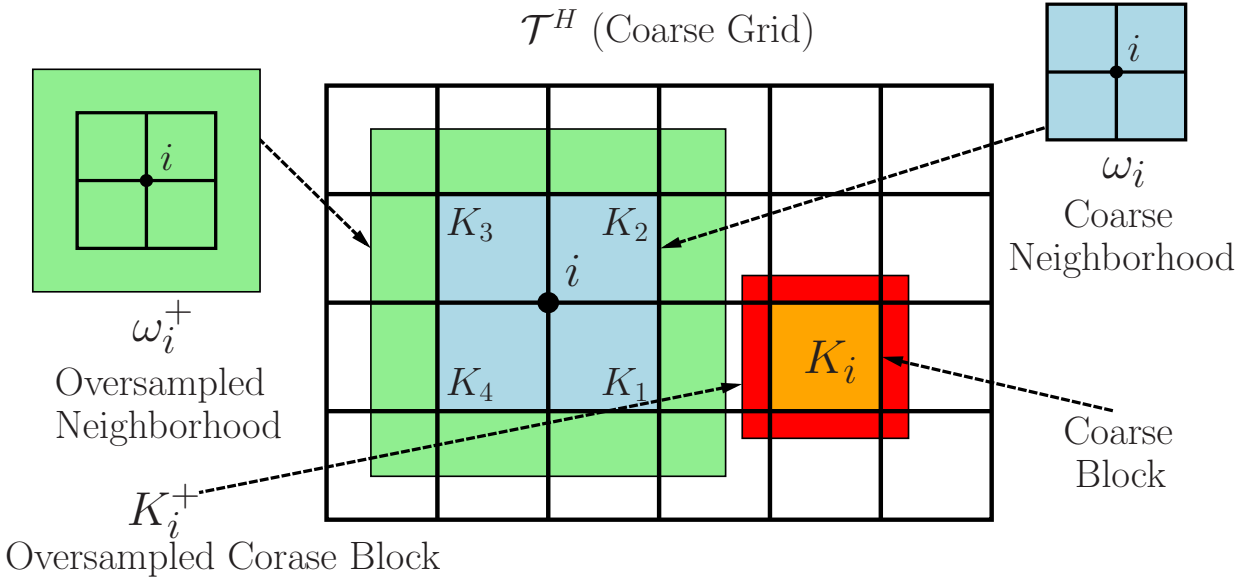}
\caption{Illustration of a coarse neighborhood, oversampled coarse neighborhood, coarse block and oversampled coarse block.}
\label{fig:grid}
\end{figure}

For each vertex $ {x}_i \in \mathcal{S}^H$ in the coarse grid, we define 
the coarse neighborhood $\omega_i$ by
\begin{equation*}
\omega_i = \bigcup \{ K_j \; : \; K_j \subset \mathcal{T}^H, \;  {x}_i \in K_j \}.
\end{equation*}
That is, $\omega_i$ is the union of all coarse grid blocks $K_j$
having the vertex $ {x}_i$
(see Figure~\ref{fig:grid}). 
A snapshot space $V^{i,\text{snap}}$ is constructed for each coarse neighborhood $\omega_i$.
The snapshot space contains a large set that represents the local solution space.
A spectral problem is then constructed to get a reduced dimensional space. 
Specifically, the spectral problem is solved in the snapshot space
and eigenfunctions corresponding to dominant modes are used 
as the final basis functions. 
To obtain conforming basis functions, each of these selected modes
will be multiplied by a partition of unity function.
The resulting space is denoted by $V^{i,\text{off}}$,
which is called the offline space for the $i$-th 
coarse neighborhood $\omega_i$.
The global offline space $V^{\text{off}}$
is then defined as the linear span of all these $V^{i,\text{off}}$,
for $i=1,2,\cdots, N_S$.
The CG coupling can be formulated as to find $ {u}_H^{\text{CG}} \in V^{\text{off}}$
such that
\begin{equation}
\label{cg_ms_sol}
a( {u}_H^{\text{CG}},  {v}) = ( {f},  {v}), \quad \forall  {v}\in V^{\text{off}}.
\end{equation}

The DG coupling can be constructed in a similar fashion.
A snapshot space $V^{i,\text{snap}}$ is constructed for each coarse grid block $K_i$.
A spectral problem is then solved in the snapshot space
and eigenfunctions corresponding to dominant modes are used as the final basis functions. 
This space is called the offline space $V^{i,\text{off}}$ for the $i$-th coarse grid block.
The global offline space $V^{\text{off}}$
is then defined as the linear span of all these $V^{i,\text{off}}$,
for $i=1,2,\cdots, N$.
The DG coupling can be formulated as: find $ {u}_H^{\text{DG}} \in V^{\text{off}}$
such that
\begin{equation}
a_{\text{DG}}( {u}_H^{\text{DG}},  {v}) = ( {f},  {v}), \quad \forall  {v}\in V^{\text{off}},
\label{eq:ipdg}
\end{equation}
where the bilinear form $a_{\text{DG}}$ is defined as
\begin{equation}
a_{\text{DG}}( {u},  {v}) = a_H( {u},  {v})
- \sum_{E\in \mathcal{E}^H} \int_E \Big( \average{ {\sigma}( {u}) \,  {n}_E} \cdot \jump{ {v}} 
+ \average{ {\sigma}( {v}) \,  {n}_E} \cdot \jump{ {u}} \Big) \; ds
+ \sum_{E\in\mathcal{E}^H} \frac{\gamma}{h} \int_E \average{\lambda+2\mu} \jump{ {u}} \cdot \jump{ {v}} \; ds
\label{eq:bilinear-ipdg}
\end{equation}
with
\begin{equation}
a_H( {u},  {v}) = 
\sum_{K\in\mathcal{T}_{H}} a_H^K(u,v),
\quad
a_H^K(u,v) = 
\int_{K}  \Big( 2\mu {\epsilon}({u}): {\epsilon}({v})
+  \lambda \nabla\cdot  {u}  \nabla\cdot {v} \Big) \; d {x},
\end{equation}
where $\gamma > 0$ is a penalty parameter, $ {n}_E$ is a fixed unit normal vector defined on the coarse edge $E$
and $ {\sigma}( {u}) \,  {n}_E$
is a matrix-vector product. 
Note that, in (\ref{eq:bilinear-ipdg}), the average and the jump operators are defined 
in the classical way.
Specifically, consider an interior coarse edge $E\in\mathcal{E}^H$
and let $K^{+}$ and $K^{-}$ be the two coarse grid blocks sharing the edge $E$.
For a piecewise smooth function $G$, we define
\begin{equation*}
\average{G} = \frac{1}{2}(G^{+} + G^{-}), \quad\quad \jump{G} = G^{+} - G^{-}, \quad\quad \text{ on } \, E,
\end{equation*}
where $G^{+} = G|_{K^{+}}$ and $G^{-} = G|_{K^{-}}$
and we assume that the normal vector $ {n}_E$
is pointing from $K^{+}$ to $K^{-}$. 
For a coarse edge $E$ lying on the boundary $\partial D$, we define
\begin{equation*}
\average{G} = \jump{G} = G, \quad\quad \text{ on } \, E,
\end{equation*}
where we always assume that $ {n}_E$ is pointing outside of $D$.
For vector-valued functions, the above average and jump operators are defined component-wise. 
We note that the DG coupling (\ref{eq:ipdg})
is the classical interior penalty discontinuous Galerkin (IPDG) method
with our multiscale basis functions.

Finally, we remark that, we use the same notations $V^{i,\text{snap}}, V^{i,\text{off}}$ and $V^{\text{off}}$
to denote the local snapshot, local offline and global offline spaces
for both the CG coupling and the DG coupling to simplify notations.

 %Let $\Omega\subset\mathbb{R}^{2}$ be a bounded domain. The displacement field $  u=(u^1,u^2)$ of matrerial is described by the following linear elasticity equations.
 
%\begin{equation}
%-\nabla\cdot \sigma(  u)=f ~~ \text{in}~\Omega,  
%\end{equation}
%\begin{equation}
 % \sigma(  u)=\mathbb C:\epsilon(  u) ~~  \text{in}~ \Omega,
%\end{equation}
%where $ \sigma$ is the stress tensor, $\epsilon$ is define as
%$\epsilon=\epsilon(  u)=\frac{1}{2}(\nabla   u+\nabla   u^T)$, $\mathbb C=\mathbb C(x),x\in\Omega$ is the fourth order elasticity tensor
%Here, we only consider the isotropic case,which means $\sigma_{ij}=\lambda\delta_{ij}e_{kk}+2\mu e_{ij}$
%\begin{equation*}
% {\sigma} = 2\mu  {\epsilon} + \lambda \nabla\cdot  {u}  {I}
%\end{equation*}

\section{Construction of multiscale basis functions}
\label{sec:gmsfem}

This section is devoted to the construction of  multiscale basis functions. 

\subsection{Basis functions for CG coupling}\label{key:cg_coupling}

We begin by the construction of local snapshot spaces. 
Let $\omega_i$ be a coarse neighborhood, $i=1,2,\cdots, N_S$. 
We will define two types of local snapshot spaces. 
The first type of local snapshot space is
\begin{equation*}
V_1^{i,\text{snap}} = V^h(\omega_i),
\end{equation*}
where $V^h(\omega_i)$ is the restriction of the conforming space to $\omega_i$. 
Therefore, $V_1^{i,\text{snap}}$ contains all possible fine scale functions defined on $\omega_i$. 
The second type of local snapshot space contains all possible harmonic extensions. Next, let $V^h(\partial\omega_i)$ be the restriction of the conforming space to $\partial\omega_i$.
Then we define the fine-grid delta function $\delta_k \in V^h(\partial\omega_i)$ on $\partial\omega_i$ by
\begin{equation*}
\delta_k( {x}_l) = 
\begin{cases}
1, \quad & l = k \\
0, \quad & l \ne k,
\end{cases}
\end{equation*}
where $ \{{x}_l\}$ are all fine grid nodes on $\partial\omega_i$. 
Given $\delta_k$, we find $ {u}_{k1}$ and $ {u}_{k2}$ by
\begin{equation}
\begin{split}
- \nabla \cdot  {\sigma}( {u}_{k1}) &=  {0}, \quad\text{ in } \; \omega_i \\
 {u}_{k1} &= (\delta_k, 0)^T, \quad\text{ on } \; \partial\omega_i
\end{split}
\label{eq:cg_snap_har_1}
\end{equation}
and 
\begin{equation}
\begin{split}
- \nabla \cdot  {\sigma}( {u}_{k2}) &=  {0}, \quad\text{ in } \; \omega_i \\
 {u}_{k2} &= (0,\delta_k)^T, \quad\text{ on } \; \partial\omega_i.
\end{split}
\label{eq:cg_snap_har_2}
\end{equation}
The linear span of the above harmonic extensions is our second type of local snapshot space $V^{i,\text{snap}}_2$. 
To simplify the notations, we will use $V^{i,\text{snap}}$ to denote $V^{i,\text{snap}}_1$ or $V^{i,\text{snap}}_2$
when there is no need to distinguish the two type of spaces. 
Moreover, we write
\begin{equation*}
V^{i,\text{snap}} = \text{span} \{  {\psi}^{i,\text{snap}}_k, \quad k=1,2,\cdots, M^{i,\text{snap}} \},
\end{equation*}
where $M^{i,\text{snap}}$ is the number of basis functions in $V^{i,\text{snap}}$.

We will perform a dimension reduction on the above snapshot spaces
by the use of a spectral problem. 
First, we will need a partition of unity function $\chi_i$
for the coarse neighborhood $\omega_i$. 
One choice of a partition of unity function
is the coarse grid hat functions $\Phi_i$, that is, the piecewise bi-linear function on the coarse grid
having value $1$ at the coarse vertex $ {x}_i$
and value $0$ at all other coarse vertices. 
The other choice is the multiscale partition of unity function, which is defined in the following way.
Let $K_j$ be a coarse grid block having the vertex $ {x}_i$. Then we consider
\begin{equation}
\begin{split}
- \nabla \cdot  {\sigma}( {\zeta}_i) &=  {0}, \quad\text{ in } \; K_j \\
 {\zeta}_{i} &= (\Phi_i,0)^T, \quad\text{ on } \; \partial K_j.
\end{split}
\end{equation}
Then we define the multiscale partition of unity as $\widetilde{\Phi}_i = ( {\zeta}_i)_1$. 
The values of $\widetilde{\Phi}_i $ on the other coarse grid blocks are defined similarly. 

Based on our analysis to be presented in the next sections, 
we define the spectral problem as
\begin{equation}
\int_{\omega_i} \Big( 2\mu  {\epsilon}( {u}) :  {\epsilon}( {v}) 
+ \lambda \nabla\cdot  {u} \, \nabla\cdot  {v}  \Big) \; d {x}
= \xi \int_{\omega_i} \tilde{\kappa}   {u} \cdot  {v} \; d {x},
\label{eq:spec-cg}
\end{equation}
where $\xi$ denotes the eigenvalue and 
\begin{equation}
\tilde\kappa = \sum_{i=1}^{N_S}(\lambda+2\mu) | \nabla \chi_i |^2.
\label{eq:kappa_tilda}
\end{equation}

The above spectral problem (\ref{eq:spec-cg}) is solved in the snapshot space. 
We let $(\phi_k, \xi_k)$ 
be the eigenfunctions and the corresponding eigenvalues. 
Assume that
\begin{equation*}
\xi_1 \leq \xi_2 \leq \cdots \leq\xi_{M^{i,\text{snap}}}.
\end{equation*}
Then the first $L_i$ eigenfunctions will be used to construct the local offline space. 
We define
\begin{equation}
 {\psi}^{i,\text{off}}_l = \sum_{k=1}^{M^{i,\text{snap}}} \phi_{lk}  {\psi}^{i,\text{snap}}_k, \quad\quad l=1,2,\cdots, L_i,
\end{equation}
where $\phi_{lk}$ is the $k$-th component of $\phi_l$. 
The local offline space is then defined as
\begin{equation*}
V^{i,\text{off}} = \text{span} \{  \chi_i  {\psi}^{i,\text{off}}_l, \quad l=1,2,\cdots, L_i \}.
\end{equation*}
Next, we define the global continuous Galerkin offline space as
\begin{equation*}
V^{\text{off}} = \text{span} \{  V^{i,\text{off}}, \quad i=1,2,\cdots, N_S \}.
\end{equation*}

\subsection{Basis functions for DG coupling}\label{key:dg_coupling}

We will construct the local basis functions required for the DG coupling. 
We also provide two types of snapshot spaces as in CG case. 
The first type of local snapshot space is
all possible fine grid bi-linear functions defined on $K_i$.
The second type of  local snapshot space $V^{i,snap}$ for the coarse grid block $K_i$
is defined as the linear span of all harmonic extensions. 
Specifically, given $\delta_k$, we find $ {u}_{k1}$ and $ {u}_{k2}$ by
\begin{equation}
\begin{split}
- \nabla \cdot  {\sigma}( {u}_{k1}) &=  {0}, \quad\text{ in } \; K_i \\
 {u}_{k1} &= (\delta_k, 0)^T, \quad\text{ on } \; \partial K_i
\end{split}
\label{eq:dg_snap_har_1}
\end{equation}
and 
\begin{equation}
\begin{split}
- \nabla \cdot  {\sigma}( {u}_{k2}) &=  {0}, \quad\text{ in } \; K_i \\
 {u}_{k2} &= (0,\delta_k)^T, \quad\text{ on } \; \partial K_i.
\end{split}
\label{eq:dg_snap_har_2}
\end{equation}
The linear span of the above harmonic extensions is the local snapshot space $V^{i,\text{snap}}$. 
We also write
\begin{equation*}
V^{i,\text{snap}} = \text{span} \{  {\psi}^{i,\text{snap}}_k, \quad k=1,2,\cdots, M^{i,\text{snap}} \},
\end{equation*}
where $M^{i,\text{snap}}$ is the number of basis functions in $V^{i,\text{snap}}$.

We will perform a dimension reduction on the above snapshot spaces
by the use of a spectral problem. 
Based on our analysis to be presented in the next sections, 
we define the spectral problem as
\begin{equation}
\int_{K_i} \big(2\mu {\epsilon}( {u}): {\epsilon}( {v}) +
 \lambda \nabla\cdot  {u}  \nabla\cdot {v}\big)d {x} 
= \frac{\xi}{H} \int_{\partial K_i}\left\langle\lambda+2\mu\right\rangle  {u}\cdot  {v} \; ds,
\label{eq:spec-dg}
\end{equation}
where $\xi$ denotes the eigenvalues and $\left\langle\lambda+2\mu\right\rangle$ is the maximum value of $\average{ \lambda+2\mu}$ on $\partial K_i$.
The above spectral problem (\ref{eq:spec-dg}) is again solved in the snapshot space $V^{i,\text{snap}}$. 
We let $(\phi_k, \xi_k)$, for $k=1,2,\cdots, M^{i,\text{snap}}$
be the eigenfunctions and the corresponding eigenvalues. 
Assume that
\begin{equation*}
\xi_1\leq \xi_2 \leq \cdots \leq \xi_{M^{i,\text{snap}}}.
\end{equation*}
Then the first $L_i$ eigenfunctions will be used to construct the local offline space. 
Indeed, we define
\begin{equation}
 {\psi}^{i,\text{off}}_l = \sum_{k=1}^{M^{i,\text{snap}}} \phi_{lk}  {\psi}^{i,\text{snap}}_k, \quad\quad l=1,2,\cdots, L_i,
\end{equation}
where $\phi_{lk}$ is the $k$-th component of $\phi_l$. 
The local offline space is then defined as
\begin{equation*}
V^{i,\text{off}} = \text{span} \{   {\psi}^{i,\text{off}}_l, \quad l=1,2,\cdots, L_i \}.
\end{equation*}
The global offline space is also defined as 
\begin{equation*}
V^{\text{off}} = \text{span} \{  V^{i,\text{off}}, \quad i=1,2,\cdots, N \}.
\end{equation*}

\subsection{Oversampling technique}

In this section,
 we present an
oversampling technique for generating multiscale basis functions. The main idea
of
oversampling is to solve local spectral problem in a larger domain. 
This allows obtaining a snapshot space that has a smaller dimension since snapshot
vectors contain solution oscillation near the boundaries. In our previous approaches,
we assume that the snapshot vectors can have an arbitrary value on the boundary
of coarse blocks which yield to large dimensional coarse spaces.

For the harmonic extension snapshot case, 
we  solve equation (\ref{eq:cg_snap_har_1}) and (\ref{eq:cg_snap_har_2}) in $\omega_i^+$ (see Figure \ref{fig:grid})
instead of $\omega_i$ for CG case, and solve the equation (\ref{eq:dg_snap_har_1}) and (\ref{eq:dg_snap_har_2}) in $K_i^+$ instead of
$K_i$ for DG case. We denote the solutions as $\psi_i^{+,\text{snap}}$, and their restrictions on $\omega_i$ or $ K_i $ as $\psi_i^{\text{snap}}$. 
We reorder these functions according eigenvalue behavior and write
$$
R_{\text{snap}}^+ = \left[ \psi_{1}^{+,\text{snap}}, \ldots, \psi_{M_{\text{snap}}}^{+,\text{snap}} \right] \quad \text{and} \quad R_{\text{snap}} = \left[ \psi_{1}^{\text{snap}}, \ldots, \psi_{M_{\text{snap}}}^{\text{snap}} \right].
$$
where $M_{snap}$ denotes the total number of functions kept in the snapshot space.
%We also propose different kind of spectral problems on offline stage. 

For CG case we define the following spectral problems in the space of snapshot:

\begin{equation}
R_{\text{snap}}^TAR_{\text{snap}}\Psi_k=\zeta (R_{\text{snap}}^{+})^TM^+R_{\text{snap}}^+\Psi_k\label{cg_over_1},
\end{equation}
or \begin{equation}
(R_{\text{snap}}^{+})^TA^+R_{\text{snap}}^+\Psi_k=\zeta (R_{\text{snap}}^{+})^TM^+R_{\text{snap}}^+\Psi_k\label{cg_over_2},
\end{equation}
where 
\begin{equation*}
\begin{split}
& A=[a_{kl}] =\int_{\omega_i} \Big( 2\mu  {\epsilon}( {\psi_k^{\text{snap}}}) :  {\epsilon}( {\psi_l^{\text{snap}}})
+ \lambda \nabla\cdot  {\psi_k^{\text{snap}}} \, \nabla\cdot  {\psi_l^{\text{snap}}}  \Big) \; d {x} ,\\
& A^+=[a_{kl}^+] =\int_{\omega_i^+} \Big( 2\mu  {\epsilon}( {\psi_k^{+,\text{snap}}}) :  {\epsilon}( {\psi_l^{+,\text{snap}}})
+ \lambda \nabla\cdot  {\psi_k^{+,\text{snap}}} \, \nabla\cdot  {\psi_l^{+,\text{snap}}}  \Big) \; d {x} ,\\
& M^+=[m_{kl}^+] =\int_{\omega_i^+} \tilde{\kappa}   {\psi_k^{+,\text{snap}}} \cdot  {\psi_l^{+,\text{snap}}} \; d {x},
\end{split}
\end{equation*}
where $\tilde{\kappa}$ is defined through (\ref{eq:kappa_tilda}).

The local spectral problem for DG coupling is defined as
\begin{equation}
(R_{\text{snap}}^{+})^TA^+R_{\text{snap}}^+\Psi_k=\zeta (R_{\text{snap}}^{+})^TM_1^+R_{\text{snap}}^+\Psi_k\label{dg_over_1}
\end{equation}
or
\begin{equation}
(R_{\text{snap}}^{+})^TA^+R_{\text{snap}}^+\Psi_k=\zeta (R_{\text{snap}}^{+})^TM_2^+R_{\text{snap}}^+\Psi_k\label{dg_over_2}
\end{equation}
in the snapshot space, where 

\begin{equation*}
\begin{split}
& A^+=[a_{kl}^+] =\int_{K_i^+} \Big( 2\mu  {\epsilon}( {\psi_k^{+,\text{snap}}}) :  {\epsilon}( {\psi_l^{+,\text{snap}}}) 
+ \lambda \nabla\cdot  {\psi_k^{+,\text{snap}}} \, \nabla\cdot  {\psi_l^{+,\text{snap}}}  \Big) \; d {x} ,\\
& M_1^+=[m_{1,kl}^+] =\frac{1}{H}\int_{K_i^+}  \average{\lambda+2\mu}  {\psi_k^{+,\text{snap}}} \cdot  {\psi_l^{+,\text{snap}}} \; d {x}, \\
& M_2^+=[m_{2,kl}^+] =\frac{1}{H}\int_{\partial K_i^+} \average{\lambda+2\mu}    {\psi_k^{+,\text{snap}}} \cdot  {\psi_l^{+,\text{snap}}} \; d {x}.
\end{split}
\end{equation*}

After solving above local spectral problems, we form 
the offline space as in the no oversampling case, see Section \ref{key:cg_coupling} 
for CG coupling and Section \ref{key:dg_coupling} for DG coupling.

%{\bf Discuss oversampling strategy here with formal notations.}
%
%{\bf Add oversampling results in Sections 6.1. and 6.2. You can mix them with
%other results}

% % error estimate

%%%%% numerical result
 \section{Numerical result}
 \label{sec:gmsfem_num_res}

\begin{figure}[H]
\begin{centering}
%\centering
\includegraphics[scale=0.5]{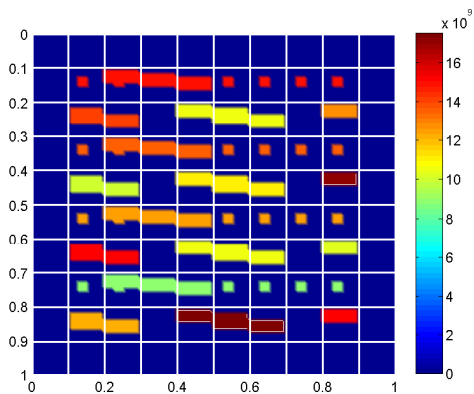}
\par\end{centering}
\caption{ $\text{Young's modulus (Model 1)}$  }
\label{fig:young_modulus}
\end{figure}

\begin{figure}[H]
\begin{centering}
%\centering
\includegraphics[scale=0.5]{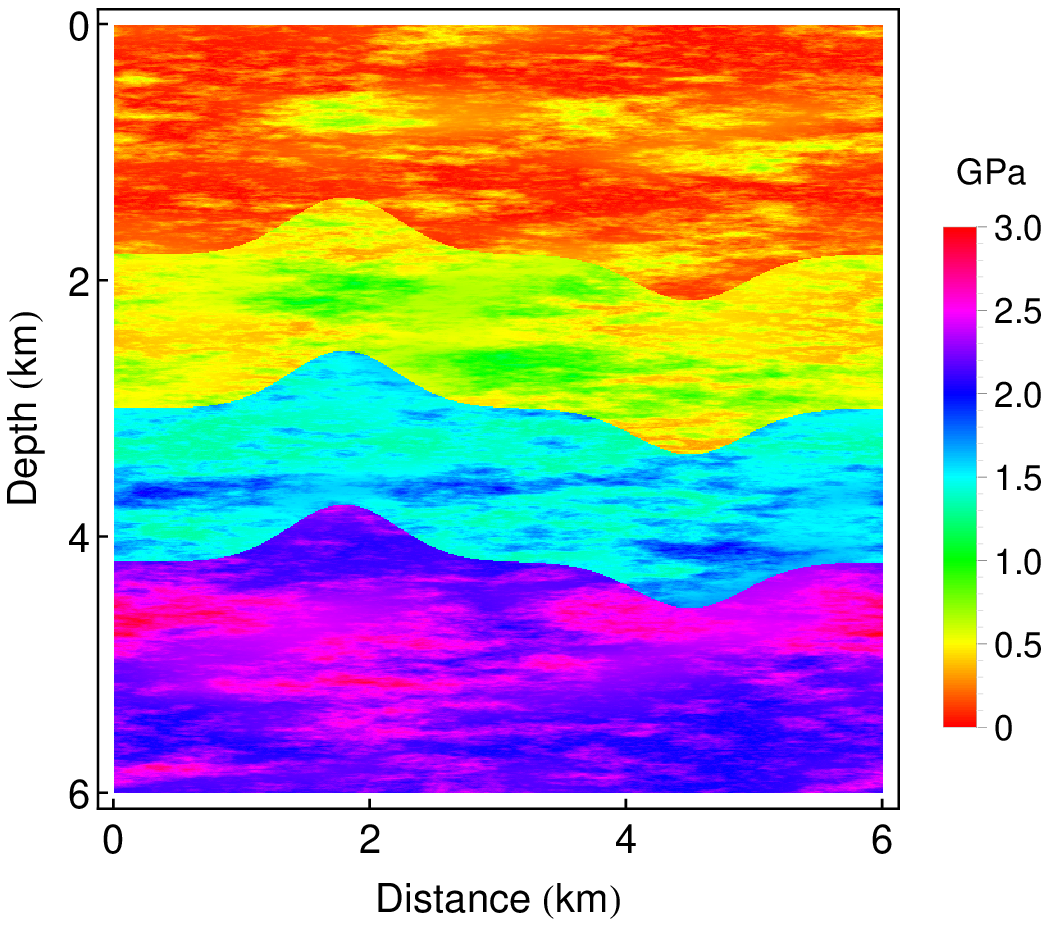}\includegraphics[scale=0.5]{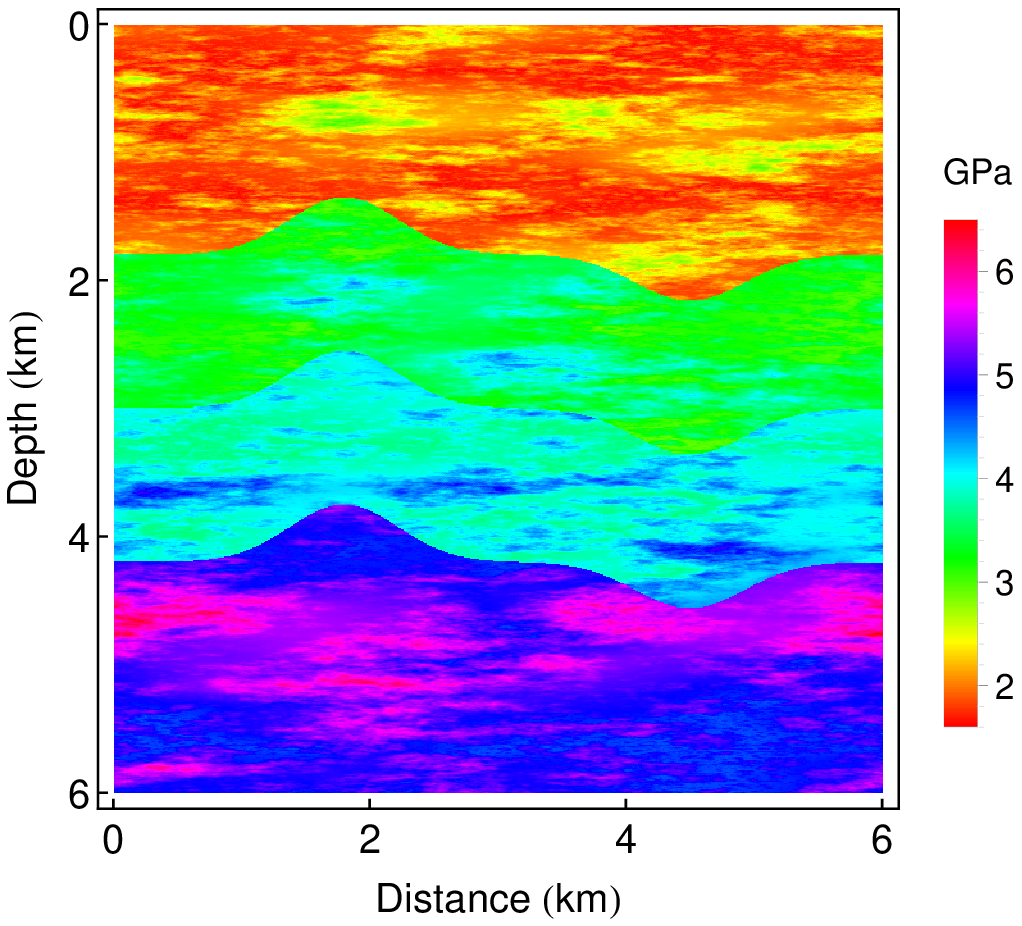}
\par\end{centering}
\caption{Left: $\lambda$ ~~Right: $\mu$ (Model 2) }
\label{fig:seg_model}
\end{figure}

In this section, we present numerical results for CG-GMsFEM and DG-GMsFEM with two models.
We consider different choices of snapshot spaces such as local-fine grid functions
and harmonic functions and use different
local spectral problems such as no-oversampling and oversampling described in the paper.
For the first model, we consider the medium that has no-scale separation and features such
as high conductivity channels and isolated inclusions.
The Young's modulus $E(x)$ is depicted in Figure \ref{fig:young_modulus}, $\lambda(x)=\frac{\nu}{(1+\nu)(1-2\nu)}E(x)$, $\mu(x)=\frac{1}{2(1+\nu)}E(x)$,
the Poisson ratio $\nu$ is taken to be $0.22$. 
For the second example, we use the model that is used in 
 \cite{gfgce14} for the simulation of subsurface elastic waves
(see Figure~\ref{fig:seg_model}).
In all numerical tests, we use constant force and homogeneous Dirchlet boundary condition. In all tables below, $\Lambda_*$  represent the minimum discarded eigenvalue of the corresponding spectral problem. We note that  the first  three eigenbasis are constant and linear functions, therefore we present our numerical results starting from fourth eigenbasis
in all cases.
%{\it In this paragraph, state the domain size, boundary conditions, and the fact that you are using heterogeneous media properties. Refer to the figure for media properties.}

%{\it In this paragraph, we will state some summary of the methods. I will write this paragraph, so skip to others}

Before presenting the numerical results, we summarize our numerical findings.

\begin{itemize}

\item We observe a fast decay in the error as more basis functions are added in both CG-GMsFEM and DG-GMsFEM

%\item We observe a better accuracy in the case of DG-GMsFEM

\item We observe the use of multiscale partition of unity improves the accuracy of
CG-GMsFEM compared to the use of piecewise bi-linear functions

\item We observe an improvement in the accuracy (a slight improvement
in CG case and a large improvement in DG case)
 when using oversampling for
the examples we considered and the decrease in the snapshot space dimension

\end{itemize}

\subsection{Numerical results for Model 1 with conforming GMsFEM (CG-GMsFEM)}

For the first model, we divide the domain $D=[0,1]\times[0,1]$ into $10\times 10$ coarse grid blocks, inside each coarse block we use $10\times 10$ fine scale square blocks, which result in a $100\times100$ fine grid blocks. 
The dimension of the reference solution is 20402. We will show the performance of CG-GMsFEM with the use of local fine-scale snapshots and harmonic extension snapshots. Both bi-linear and multiscale partition of unity functions (see section \ref{key:cg_coupling}) will be considered. For each case, we will provide the comparsion  using oversampling and no-oversampling.
%{\it In this paragraph, write down your fine mesh size, coarse mesh size. Write down that we will consider two types of partition of unity functions and refer to them. Write down that we will be using two types of snapshot spaces and state them. Write down that we will use both oversampling and no-oversampling and refer to them in appropriate Sections.}
For the error measure, we use relative weighted $L^2$ norm error and weighted $H^1$ norm error to compare the accuracy of CG-GMsFEM, which is defined as
\[
 e_{L^2}=\cfrac{\|(\lambda+2\mu )(  u_{H}-  u_h)\|_{L^{2}(D)}}{\|(\lambda+2\mu )  u_h\|_{L^{2}(D)}},\quad
 e_{H^{1}}=\sqrt{\cfrac{a( {u_{H}-u_h}, {u_{H}-u_h})}{a( {u_h}, {u_h})}}\]
where $ {u_H}$ and $ {u_h}$ are CG-GMsFEM defined in (\ref{cg_ms_sol}) and fine-scale CG-FEM solution defined in (\ref{cg_fine_sol}) respectively.

Tables \ref{lin_spec} and \ref{ms_spec} show the numerical results of using local fine-scale
snapshots with piecewise bi-linear function and multiscale functions as partition of unity respectively. As we observe, when using more multiscale basis, the errors decay rapidly, especially for multiscale partition of unity. For example, we can see that the weighted $L^2$ error drops from 24.9\% to 1.1\% in the case of using bi-linear function as partition of unity with no oversampling, while the dimension increases from 728 to 2672. If we use multiscale partition of unity, the corresponding weighted $L^2$ error drops from 8.4\% to 0.6\%, which demonstrates a great advantage of multiscale partition of unity. Oversampling can help improve the accuracy as our results indicate. The local eigenvalue problem used for oversampling is Eq.(\ref{cg_over_2}).
%\marginpar{We need to write what oversampling is used. What is the dimension of the snapshot space. Add a comment on a periodic case. In periodic case, we observe a substantial reduction in the snapshot space dimension. We have considered. Obtain the same accuracy. }

%{\it In this paragraph, describe numerical results for snapshot space 1, which are fine-scale vectors and describe the results in Tables 1-2. Emphasize that the error decreases as we increase the number of basis functions. Mention that oversampling helps. State some numbers, e.g., the error decreases to 10 percent if we use 16 basis.}

%{\it In this paragraph, describe the results for harmonic snapshot space. You can be short. Just emphasize the error decreases. And mention any numbers or construction parameters that you have}

Next, we present the numerical results when harmonic extensions are used  as snapshots in Tables \ref{lin_har} and \ref{ms_har}. We can observe similar trends as in the local fine-scale snapshot case. The errors  decrease as the number of basis functions increase. The $L^2$ error is less than $1$\%  when about $13$\% percent of degrees of freedom is used. Similarly, the oversampling method helps to 
improve the accuracy. In this case, the local eigenvalue problem used for oversampling is Eq.(\ref{cg_over_1}).

\begin{table}
\centering
\begin{tabular}{|c|c|c|c|c|c|c|}
\hline 
\multirow{2}{*}{\specialcell{Dimension}} & \multicolumn{2}{c|}{$1/\Lambda_*$} & \multicolumn{2}{c|}{$e_{L^{2}}$} & \multicolumn{2}{c|}{$e_{H^{1}}$}\tabularnewline
\cline{2-7} 
 & \specialcell{without\\oversampling} & \specialcell{with\\ oversampling} & \specialcell{without\\oversampling}  & \specialcell{with\\oversampling} & \specialcell{without\\oversampling}  & \specialcell{with\\oversampling} \tabularnewline
\hline
\hline
728  & 1.3e+07 & 1.4e+07  & 0.249 & 0.215 & 0.444 & 0.409   \tabularnewline
\hline
1214  & 3.1e+06  & 5.6e+06  & 0.048  & 0.047 & 0.220 & 0.213   \tabularnewline
\hline
1700 & 7.0e+05 & 2.7e+06  & 0.027  & 0.024& 0.162 & 0.153  \tabularnewline
\hline
2186  & 1.8e+00 & 1.7e+06  & 0.018  & 0.016  & 0.133 & 0.123   \tabularnewline
\hline
2672  & 9.9e-01 & 1.4e+06 & 0.011 & 0.010  & 0.105 & 0.099   \tabularnewline
\hline 
\end{tabular}
\caption{Relative errors between CG-MsFEM solution and the fine-scale CG-FEM solution, piecewise bi-linear partition of unity functions are used. The case with local fine-scale snapshots.}
\label{lin_spec}
\end{table}

\begin{table}
\centering
\begin{tabular}{|c|c|c|c|c|c|c|}
\hline 
\multirow{2}{*}{\specialcell{Dimension}} & \multicolumn{2}{c|}{$1/\Lambda_*$} & \multicolumn{2}{c|}{$e_{L^{2}}$} & \multicolumn{2}{c|}{$e_{H^{1}}$}\tabularnewline
\cline{2-7} 
 & \specialcell{without\\oversampling} & \specialcell{with\\oversampling} & \specialcell{without\\oversampling}  & \specialcell{with\\oversampling} & \specialcell{without\\oversampling}  & \specialcell{with\\oversampling} \tabularnewline
\hline
\hline
728   & 6.9e+06 & 6.2e+06 & 0.084& 0.110 & 0.254& 0.274  \tabularnewline
\hline
1214   & 5.8e+00& 3.2e+06 & 0.031& 0.028 & 0.166& 0.160   \tabularnewline
\hline
1700 & 2.1e+00 & 1.2e+06  & 0.015& 0.012& 0.111& 0.105  \tabularnewline
\hline
2186 & 1.3e+00 & 5.9e+05 & 0.009& 0.008  & 0.088& 0.083   \tabularnewline
\hline
2672  & 9.4e-01& 1.0e+01  & 0.006& 0.005  & 0.071& 0.066   \tabularnewline
\hline
\end{tabular}\caption{Relative errors between CG-MsFEM solution and the fine-scale CG-FEM solution, multiscale partition of unity functions are used. The case with local fine-scale snapshots.}
\label{ms_spec}
\end{table}

\begin{table}
\centering
\begin{tabular}{|c|c|c|c|c|c|c|}
\hline 
\multirow{2}{*}{\specialcell{Dimension}} & \multicolumn{2}{c|}{$1/\Lambda_*$} & \multicolumn{2}{c|}{$e_{L^{2}}$} & \multicolumn{2}{c|}{$e_{H^{1}}$}\tabularnewline
\cline{2-7} 
 & \specialcell{without\\oversampling} & \specialcell{with\\oversampling} & \specialcell{without\\oversampling}  & \specialcell{with\\oversampling} & \specialcell{without\\oversampling}  & \specialcell{with\\oversampling} \tabularnewline
\hline
\hline
728  & 1.3e+07& 1.2e+07  & 0.254& 0.218 & 0.446& 0.418  \tabularnewline
\hline
1214 & 2.1e+06& 5.5e+06  & 0.047& 0.048 & 0.218& 0.217   \tabularnewline
\hline
1700  & 2.8e+05& 3.2e+06  & 0.024& 0.022& 0.153& 0.148  \tabularnewline
\hline
2186  & 1.2e+00& 9.8e+05  & 0.016& 0.015  & 0.124& 0.122  \tabularnewline
\hline
2672 & 5.8e-01& 2.1e+04  & 0.008& 0.010  & 0.102& 0.099   \tabularnewline
\hline
\end{tabular}\caption{Relative errors between CG-MsFEM solution and the fine-scale CG-FEM solution, piecewise bi-linear partition of unity functions are used. The case with hamonic snapshots.}
\label{lin_har}
\end{table}

\begin{table}
\centering
\begin{tabular}{|c|c|c|c|c|c|c|}
\hline 
\multirow{2}{*}{Dimension} & \multicolumn{2}{c|}{$1/\Lambda_*$} & \multicolumn{2}{c|}{$e_{L^{2}}$} & \multicolumn{2}{c|}{$e_{H^{1}}$}\tabularnewline
\cline{2-7} 
 & \specialcell{without\\oversampling} & \specialcell{with\\oversampling} & \specialcell{without\\oversampling}  & \specialcell{with\\oversampling} & \specialcell{without\\oversampling}  & \specialcell{with\\oversampling} \tabularnewline
\hline
\hline
728  & 7.0e+06& 7.2e+06  & 0.087& 0.112 & 0.259& 0.291  \tabularnewline
\hline
1214 & 5.5e+00& 3.2e+06   & 0.034& 0.032 & 0.174& 0.169  \tabularnewline
\hline
1700 &  1.9e+00& 1.5e+06 & 0.015& 0.013& 0.115& 0.112  \tabularnewline
\hline
2186  &  1.0e+00& 2.5e+05 & 0.009& 0.008 & 0.090& 0.089   \tabularnewline
\hline
2672  & 7.1e-01& 1.7e+00  & 0.007& 0.006  & 0.075& 0.074   \tabularnewline
\hline
\end{tabular}\caption{Relative errors between CG-MsFEM solution and the fine-scale CG-FEM solution, multiscale partition of unity functions are used. The case with hamonic snapshots.}
\label{ms_har}
\end{table}

\subsection{Numerical results for Model 1 with DG-GMsFEM}

%{\it Follow, CG-GMsFEM Discussions and write 3-4 paragraphs.}
In this section, we consider numerical results for DG-GMsFEM discussed in Section \ref{key:dg_coupling}.
To show the performance of DG-GMsFEM, we use the same model (see Figure \ref{fig:young_modulus}) and the coarse and fine grid settings as in the CG case. We will also present the result of using both harmonic extension and eigenbasis (local fine-scale) as snapshot space. To measure the error,
 we define broken weighted $L^2$ norm error and $H^1$ norm error
\[
e_{L^2}=\sqrt{\frac{\sum_{K\in\mathcal{T}_{H}}\int_{K}(\lambda+2\mu)(  u_{H}-  u_h)^2dx}{\sum_{K\in\mathcal{T}_{H}}\int_{K} (\lambda+2\mu) u_h^2dx}} \quad
%e_{\text{en}}= \sqrt{\frac{a_h^{DG}(  u_{H}-  u,  u_{H}-  u)}{a_h^{DG}( {u}, {u_h})}}
%\]
%\[
e_{H^1} =\sqrt{\frac{\sum_{K\in\mathcal{T}_{H}}\int_{K} {\sigma}(  u_{H}-  u_h)): {\varepsilon}(  u_{H}-  u_h))dx}{\sum_{K\in\mathcal{T}_{H}}\int_{K} {\sigma}(  u_h): {\varepsilon}( u_h)dx}}
\]
where $ {u_H}$ and $ {u_h}$ are DG-GMsFEM defined in (\ref{eq:ipdg}) and fine-scale DG-FEM solution defined in
(\ref{eq:ipdgfine}) respectively.

In Table \ref{dg-spe}, the numerical results of DG-MsFEM with local fine-scale functions as the snapshot space is shown. We observe that DG-MsFEM shows a better approximation compared to CG-MsFEM if oversampling is used. 
The error decreases more rapidly as we add basis.  More specifically, the relative broken $L^2$  error and $H^1$ error decrease from $14.1$\%,  $52.5$\% to $0.2$\% and $5.8$\% respectively, while 
the degrees of freedom of the coarse system increase from $728$ to $2696$, where the latter is only $13.2$\% of the reference solution. The local eigenvalue problem used for oversampling is Eq.(\ref{dg_over_1}).

Table \ref{dg-har} shows the corresponding results when harmonic functions are used to construct the snapshot space. We observe similar errors decay trend as local fine-scale snapshots are used. Oversampling can help improve the results significantly. Although the error is very large when the dimension of coarse system is $728$ ($4$ multiscale basis is used), the error becomes very small when the dimension reaches $1728$ ($9$ multiscale basis is used). The local eigenvalue problem used for oversampling here is Eq.(\ref{dg_over_2}). We remark that oversampling can not only help decrease the error, but also decrease the dimension of the snapshot space greatly in peridoic case. 

\begin{table}
\centering
\begin{tabular}{|c|c|c|c|c|c|c|}
\hline 
\multirow{2}{*}{Dimension} & \multicolumn{2}{c|}{$1/\Lambda_*$} & \multicolumn{2}{c|}{$e_{L^{2}}$} & \multicolumn{2}{c|}{$e_{H^{1}}$}\tabularnewline
\cline{2-7} 
 & \specialcell{without\\oversampling} & \specialcell{with\\oversampling} & \specialcell{without\\oversampling}  & \specialcell{with\\oversampling} & \specialcell{without\\oversampling}  & \specialcell{with\\oversampling} \tabularnewline
\hline
\hline
728  & 4.9e-03& 1.5e-03  & 0.281& 0.141 & 0.554& 0.525  \tabularnewline
\hline
1184  & 3.0e-03& 8.5e-04 & 0.118& 0.019 & 0.439& 0.209 \tabularnewline
\hline
1728  & 2.1e-03& 5.6e-04 & 0.108& 0.012& 0.394& 0.145  \tabularnewline
\hline
2184& 1.2e-03& 3.5e-04  &  0.073& 0.007& 0.348& 0.096  \tabularnewline
\hline
2696 & 1.0e-03& 2.7e-04 &  0.056& 0.002  & 0.300& 0.058    \tabularnewline
\hline

\end{tabular}\caption{Relative errors between DG-MsFEM solution and the fine-scale DG-FEM solution. The case with local fine-scale snapshots. }
\label{dg-spe}
\end{table}

\begin{table}
\centering
\begin{tabular}{|c|c|c|c|c|c|c|}
\hline 
\multirow{2}{*}{Dimension} & \multicolumn{2}{c|}{$1/\Lambda_*$} & \multicolumn{2}{c|}{$e_{L^{2}}$} & \multicolumn{2}{c|}{$e_{H^{1}}$}\tabularnewline
\cline{2-7} 
 & \specialcell{without\\oversampling} & \specialcell{with\\oversampling} & \specialcell{without\\oversampling}  & \specialcell{with\\oversampling} & \specialcell{without\\oversampling}  & \specialcell{with\\oversampling} \tabularnewline
\hline
\hline
728  & 2.9e-01& 1.6e-01   & 0.285& 0.149 & 0.557& 0.528  \tabularnewline
\hline
1184  & 1.6e-01&  6.5e-02 & 0.193& 0.076 & 0.515& 0.366\tabularnewline
\hline
1728  & 1.0e-01&  5.4e-02 & 0.114& 0.009& 0.432& 0.155  \tabularnewline
\hline
2184& 7.1e-02 & 3.9e-02 & 0.081& 0.004& 0.326& 0.078 \tabularnewline
\hline
2696 & 6.3e-02 & 2.8e-02& 0.043& 0.002  & 0.231& 0.060   \tabularnewline
\hline
\end{tabular}\caption{Relative errors between DG-MsFEM solution and the fine-scale DG-FEM solution. The case with hamonic snapshots. }
\label{dg-har}
\end{table}
\subsection{Numerical results for Model  2}

The purpose of this example is to test a method for an earth model that
is used in \cite{gfgce14}.
The domain for the second model is $D=(0,6000)^2$ (in meters) which is divided into $900=30\times30$ square coarse grid blocks, inside each coarse block we generate $20\times 20$ fine scale square blocks. The reference solution is computed through standard CG-FEM on the resulting $600\times 600$ fine grid. We note that the dimension of the reference solution is 722402. The numerical results for CG-MsFEM and DG-MsFEM are presented in Table \ref{seg_lin_eigen} and \ref{dg_seg_eigen}
respectively. We observe the relatively low errors compared to the high contrast case and the error decrease with the dimension increase of the offline space. Both coupling methods (CG and DG) show very good approximation ability.
\begin{table}[H]
 \centering \begin{tabular}{|c|c|c|c|c|c|c|}
\hline
dimension   & $\frac{1}{\Lambda_*}$    & $e_{L^2}$  & $e_{H^{1}}$ \tabularnewline
\hline
6968  & 4.9e+00  & 3.1e-03 & 5.4e-02  \tabularnewline
\hline
8650 & 4.5e+00   & 2.7e-03 & 5.2e-02  \tabularnewline
\hline
10332 & 3.9e+00 & 2.5e-03& 4.9e-02  \tabularnewline
\hline
12014 &  3.6e+00 & 2.2e-03 & 4.7e-02   \tabularnewline
\hline

\end{tabular}\caption{Relative errors between CG-MsFEM solution and the fine-scale CG-FEM solution,  piecewise bi-linear partition of unity functions are used. The case with  local fine-scale snapshots.}
\label{seg_lin_eigen}
\end{table}

\begin{table}[H]
 \centering \begin{tabular}{|c|c|c|c|c|c|c|}
\hline
dimension   & $\frac{1}{\Lambda_*}$   & $e_{L^2}$  & $e_{H^{1}}$ \tabularnewline
\hline
7200  & 6.3e-06 & 4.1e-03 & 7.1e-02  \tabularnewline
\hline
9000 & 6.0e-06   & 4.0e-03 & 6.6e-02  \tabularnewline
\hline
10800 & 4.6e-06 & 3.8e-03& 6.3e-02  \tabularnewline
\hline
12600 &  4.5e-06 & 3.1e-03 & 5.9e-02   \tabularnewline
\hline

\end{tabular}\caption{Relative errors between DG-MsFEM solution and the fine-scale DG-FEM solution. The case with local fine-scale snapshots.}
\label{dg_seg_eigen}
\end{table}

\section{Error estimate for CG coupling}
\label{sec:gmsfem_error}

In this section, we present error analysis for both no oversampling
and oversampling cases. In the below, $a\preceq b$ means $a\leq Cb$, where $C$ is a contant that is independend of the mesh size and the contrast of the coefficient.

\subsection{No oversampling case}

\begin{lemma}
\label{cacc_lemma}
Let $\omega_n$ coarse neighborhood. For any $\psi\in H^1(\omega_n) $, we define $r=-div(\sigma(\psi))$. Then we have
\begin{equation}
\label{cacc}
\int_{\omega_n}2\mu \chi_n^2 \epsilon(\psi):\epsilon(\psi)+\int_{\omega_n}\lambda \chi_n^2 (\nabla\cdot{\psi})^2\\
\preceq |\int_{\omega_n} \chi_n^2{r}\cdot{{\psi}} |+\int_{\omega_n} (\lambda + 2\mu)|\nabla \chi_n|^2 {\psi}^2,
\end{equation}
 where $\chi_n$ is a scalar partition of unity subordinated to the coarse neighborhood $\omega_n$. %\marginpar{What is $V^h$}
\end{lemma}
{\it Proof}.
Multiplying both sides of $-div(\sigma(\psi))=r$ by $\chi_n^2 \psi$, we have
\begin{equation}
\begin{split}
\int_{\omega_n} \chi_n^2{r}\cdot{\psi} &=\int_{\omega_n} 2 \mu \epsilon(\psi):\epsilon(\chi_n^2\psi)+\int_{\omega_n}\lambda \nabla\cdot\psi \nabla\cdot(\chi_n^2{\psi})  \\
&= \int_{\omega_n} 2\mu \chi_n^2 \epsilon(\psi):\epsilon(\psi)+ \int_{\omega_n}2 \mu \chi_n \epsilon_{ij}({\psi}) (  \psi_i {\partial \chi_n\over \partial x_j} +  \psi_j {\partial \chi_n\over \partial x_i})\\
&\quad+\int_{\omega_n} \lambda \chi_n^2 (\nabla\cdot{\psi})^2 + \int_{\omega_n}2 \lambda \nabla\cdot\psi \chi_n \psi\cdot\nabla\chi_n \\
&=\int_{\omega_n} 2\mu \chi_n^2 \epsilon(\psi):\epsilon(\psi) +\int_{\omega_n} 2 \left( \sqrt{2\mu} \chi_n \epsilon_{ij}({\psi})\right) \left(\sqrt{\mu/2} (  \psi_i {\partial \chi_n\over \partial x_j} +  \psi_j {\partial \chi_n\over \partial x_i})\right)\\
&\quad+\int_{\omega_n}\lambda \chi_n^2 (\nabla\cdot\psi)^2 +
\int_{\omega_n} 2 \left( \sqrt{\lambda} \chi_n \nabla\cdot\psi\right) \left(\sqrt{\lambda}\psi\cdot\nabla\chi_n\right).
\end{split}
\end{equation}
Therefore,
\begin{equation}
\begin{split}\label{caco}
& \int_{\omega_n}2\mu \chi_n^2 \epsilon(\psi):\epsilon(\psi) +\int_{\omega_n}\lambda \chi_n^2 (\nabla\cdot\psi)^2\\
&\leq|\int_{\omega_n} \chi_n^2{r}\cdot{\psi} |+|
\int_{\omega_n} 2 \left( \sqrt{2\mu} \chi_n \epsilon_{ij}({\psi})\right) \left(\sqrt{\mu/2} (  \psi_i {\partial \chi_n\over \partial x_j} +  \psi_j {\partial \chi_n\over \partial x_i})\right)+
\int_{\omega_n} 2 \left( \sqrt{\lambda} \chi_n \nabla\cdot\psi\right) \left(\sqrt{\lambda}\psi\cdot\nabla\chi_n\right)|\\
&\preceq |\int_{\omega_n} \chi_n^2{r}\cdot{{\psi}} |+\int_{\omega_n} (2\lambda + 4\mu)|\nabla \chi_n|^2 {\psi}^2\\
&\preceq |\int_{\omega_n} \chi_n^2{r}\cdot{{\psi}} |+\int_{\omega_n} (\lambda + 2\mu)|\nabla \chi_n|^2 {\psi}^2.
\end{split}
\end{equation}
In the last step, we have used $2ab\leq \epsilon a^2+\frac{1}{\epsilon}b^2$, and $(ab+cd)^2\leq (a^2+c^2)(b^2+d^2)$.

\begin{flushright}
$\square$
\end{flushright}

Next, we will show the convergence of the CG-GMsFEM solution defined in (\ref{cg_ms_sol}) without oversampling. We take $I^{\omega_n}{u_h}$ to be the first $L_n$ terms of spectral expansion of $u$ in terms of eigenfunctions of the problem
$ -\text{div}(\sigma(\phi_{n}))=\xi \tilde{\kappa}\phi_{n}$ solved in $V^h{(\omega_n)}$. Applying Cea's Lemma, Lemma \ref{cacc_lemma} and using the fact that $\chi_n\preceq 1$, we can get
\begin{equation}
\label{cg_equ1}
\begin{split}
& \quad\int_D \left(2\mu \epsilon({u_h}-{u_H}):\epsilon({u_h}-{u_H})+\lambda (\nabla\cdot({u_h}-{u_H}))^2\right)  \\
&\preceq \sum_{n=1}^{N_s} \int_{\omega_n}\left(2\mu  \epsilon(\chi_n({u_h}-I^{\omega_n}{u_h})):\epsilon(\chi_n({u_h}-I^{\omega_n}{u_h})) +\lambda  (\nabla\cdot(\chi_n({u_h}-I^{\omega_n}{u_h})))^2\right) \\ 
&\preceq\sum_{n=1}^{N_s}\int_{\omega_n}2\mu \chi_n^2 \epsilon({u_h}-I^{\omega_n}{u_h}):\epsilon({u_h}-I^{\omega_n}{u_h}) +\sum_{n=1}^{N_s}\int_{\omega_n}\lambda \chi_n^2 (\nabla\cdot({u_h}-I^{\omega_n}{u_h}))^2\\
&\quad+\sum_{n=1}^{N_s} \int_{\omega_n}(\lambda +2\mu)|\nabla \chi_n|^2 ({u_h}-I^{\omega_n}{u_h})^2\\
&\preceq\sum_{n=1}^{N_s} \int_{\omega_n}(\lambda +2\mu)|\nabla \chi_n|^2 ({u_h}-I^{\omega_n}{u_h})^2+\sum_{n=1}^{N_s}|\int_{\omega_n} \chi_n^2g\cdot({u_h}-I^{\omega_n}{u_h}) |\\
&\preceq\sum_{n=1}^{N_s} \int_{\omega_n}(\lambda +2\mu)|\nabla \chi_n|^2 ({u_h}-I^{\omega_n}{u_h})^2+\sum_{n=1}^{N_s}\int_{\omega_n}((\lambda+2\mu) |\nabla\chi_n|^2)^{-1} {g}^2,
\end{split}
\end{equation}
where $g=f+\text{div}(\sigma(I^{\omega_n}{u_h}))$, $f$ is the right hand side of (\ref{ob_equ}).%$g$ is the residual in the approximation.\\

Using the properties of the eigenfunctions, we obtain
\begin{equation}
\int_{\omega_n}(\lambda+2\mu)\sum_{s=1}^{N_s}|\nabla\chi_s|^2 ({u_h}-I^{\omega_n}{u_h})^2\preceq \frac{1}{\xi_{L_n+1}^{\omega_n}}\left(\int_{\omega_n}2\mu  \epsilon({u_h}-I^{\omega_n}{u_h}):\epsilon({u_h}-I^{\omega_n}{u_h})+ 
\lambda  (\nabla\cdot({u_h}-I^{\omega_n}{u_h}))^2\right)  . 
\end{equation}
Then, the first term in the right hand side of (\ref{cg_equ1}) can be estimated as follows
\begin{equation}
\label{cg_eq2}
\begin{split}
&\quad\sum_{n=1}^{N_s} \int_{\omega_n}(\lambda +2\mu)|\nabla \chi_n|^2 ({u_h}-I^{\omega_n}{u_h})^2
\preceq \sum_{n=1}^{N_s}\int_{\omega_n}(\lambda+2\mu)\sum_{s=1}^{N_s}|\nabla\chi_s|^2 |({u_h}-I^{\omega_n}{u_h})^2\\
&\preceq \sum_{n=1}^{N_s}\frac{1}{\xi_{L_n+1}^{\omega_n}}\left(\int_{\omega_n}2\mu  \epsilon({u_h}-I^{\omega_n}{u_h}): \epsilon({u_h}-I^{\omega_n}{u_h})+\int_{\omega_n}\lambda  (\nabla\cdot({u_h}-I^{\omega_n}{u_h}))^2\right) \\
&\preceq \sum_{n=1}^{N_s}\frac{\alpha^{\omega_n}_{L_n+1}}{\xi_{L_n+1}^{\omega_n}}\left(\int_{\omega_n}2\mu \chi_n^2 \epsilon({u_h}-I^{\omega_n}{u_h}):\epsilon({u_h}-I^{\omega_n}{u_h})+\int_{\omega_n}\lambda \chi_n^2 (\nabla\cdot({u_h}-I^{\omega_n}{u_h}))^2\right)   \\
&\preceq \sum_{n=1}^{N_s}\frac{\alpha^{\omega_n}_{L_n+1}}{\xi_{L_n+1}^{\omega_n}}\int_{\omega_n} (\lambda + 2\mu)|\nabla \chi_n|^2 ({u_h}-I^{\omega_n}{u_h})^2+
\sum_{n=1}^{N_s}\frac{\alpha^{\omega_n}_{L_n+1}}{\xi_{L_n+1}^{\omega_n}}|\int_{\omega_n} \chi_n^2{g}\cdot({u_h}-I^{\omega_n}{u_h}) |\\
&\preceq \frac{1}{\Lambda_*}\left(\sum_{n=1}^{N_s}\int_{\omega_n}(\lambda + 2\mu)|\nabla \chi_n|^2 ({u_h}-I^{\omega_n}{u_h})^2+\sum_{n=1}^{N_s}|\int_{\omega_n} \chi_n^2{g}\cdot({u_h}-I^{\omega_n}{u_h}) |\right),
\end{split}
\end{equation}
where
\begin{equation*}
\Lambda_*=\text{min}_{\omega_n}\frac{\xi_{L_n+1}^{\omega_n}}{\alpha^{\omega_n}_{L_n+1}},
\end{equation*}
and\\
\begin{equation*}
{\alpha^{\omega_n}_{L_n+1}}=\frac{\int_{\omega_n} 2\mu  \epsilon({u_h}-I^{\omega_n}{u_h}):\epsilon({u_h}-I^{\omega_n}{u_h})+\int_{\omega_n}\lambda (\nabla\cdot({u_h}-I^{\omega_n}{u_h}))^2} {\int_{\omega_n}2\mu \chi_n^2 \epsilon({u_h}-I^{\omega_n}{u_h}):\epsilon({u_h}-I^{\omega_n}{u_h})+\int_{\omega_n}\lambda \chi_n^2 (\nabla\cdot({u_h}-I^{\omega_n}{u_h}))^2}.
\end{equation*}

Applying inequality (\ref{cg_eq2}) $m$ times, we have
\begin{equation}
\begin{split}
&\quad\sum_{n=1}^{N_s} \int_{\omega_n}(\lambda +2\mu)|\nabla \chi_n|^2 ({u_h}-I^{\omega_n}{u_h})^2\\
&\preceq \left(\frac{1}{\Lambda_*}\right)^m\sum_{n=1}^{N_s}\int_{\omega_n}(\lambda + 2\mu)|\nabla \chi_n|^2 ({u_h}-I^{\omega_n}{u_h})^2+\sum_{l=1}^m \left(\frac{1}{\Lambda_*}\right)^l
\sum_{n=1}^{N_s}|\int_{\omega_n} \chi_n^2{g}\cdot({u_h}-I^{\omega_n}{u_h}) |\\
%%&\preceq \left(\frac{1}{\Lambda_*}\right)^m\sum_{n=1}^{N_s}\int_{\omega_n}(\lambda + 2\mu)|\nabla \chi_n|^2 ({u_h}-I^{\omega_n}{u_h})^2+\left(\frac{1-\Lambda_*^{-m}}{\Lambda_*-1}\right)
%\sum_{n=1}^{N_s}|\int_{\omega_n} \chi_n^2{g}\cdot({u_h}-I^{\omega_n}{u_h}) |\\
&\preceq \left(\frac{1}{\Lambda_*}\right)^m\sum_{n=1}^{N_s}\int_{\omega_n}(\lambda + 2\mu)|\nabla \chi_n|^2 ({u_h}-I^{\omega_n}{u_h})^2+ (\Lambda_*)^m\left(\frac{1-\Lambda_*^{-m}}{\Lambda_*-1}\right)\sum_{n=1}^{N_s}\int_{\omega_n}((\lambda+2\mu)|\nabla\chi_n|^2)^{-1}{g}^2,
\end{split}
\end{equation}
Taking into account that
\begin{equation}
\sum_{n=1}^{N_s} \int_{\omega_n}(\lambda +2\mu)|\nabla \chi_n|^2 ({u_h}-I^{\omega_n}{u_h})^2\preceq \sum_{n=1}^{N_s}\int_{\omega_n}(\lambda+2\mu)\sum_{s=1}^{N_s} |\nabla\chi_s|^2 ({u_h}-I^{\omega_n}{u_h})^2,\\
\end{equation}
and
\begin{equation}
\label{fem_estimate}
\sum_{n=1}^{N_s}\int_{\omega_n} \left(2\mu  \epsilon({u_h}-I^{\omega_n}{u_h}):\epsilon({u_h}-I^{\omega_n}{u_h})+\lambda  (\nabla\cdot({u_h}-I^{\omega_n}{u_h}))^2\right) \preceq \int_{D}\left(2\mu \epsilon({u_h}):\epsilon(u)+\lambda (\nabla\cdot{u_h})^2\right) .
\end{equation}
inequality (\ref{cg_equ1}) becomes 
\begin{equation}
\begin{split}
&\quad\int_D\left(2\mu  \epsilon({u_h}-{u_H}): \epsilon({u_h}-{u_H})+\lambda (\nabla\cdot({u_h}-{u_H}))^2\right)  \\
& \preceq\left(\frac{1}{\Lambda_*}\right)^{m+1}\left(\sum_{n=1}^{N_s}\int_{\omega_n}2\mu  \epsilon({u_h}-I^{\omega_n}{u_h}):\epsilon({u_h}-I^{\omega_n}{u_h}) +\sum_{n=1}^{N_s}\int_{\omega_n}\lambda  (\nabla\cdot({u_h}-I^{\omega_n}{u_h}))^2\right)  \\
&\quad+\left(\Lambda_*^m\left(\frac{1-\Lambda_*^{-m}}{\Lambda_*-1}\right)+1\right)
\sum_{n=1}^{N_s}\int_{\omega_n}((\lambda+2\mu)|\nabla\chi_n|^2)^{-1}{g}^2\\
&\preceq\big(\frac{1}{\Lambda_*}\big)^{m+1}\int_{D}  \left(2\mu  \epsilon({u_h}):\epsilon(u)+\lambda(\nabla\cdot{u_h})^2\right)
+ \left((\Lambda_*)^m\left(\frac{1-(\Lambda_*)^{-m}}{\Lambda_*-1}\right)+1\right)R,
\end{split}
\end{equation}
where $R=\sum_{n=1}^{N_s}\int_{\omega_n}((\lambda+2\mu)|\nabla\chi_n|^2)^{-1}{g}^2$. If $|{g}|\preceq1$, then $\int_{\omega_n}((\lambda+2\mu)|\nabla\chi_n|^2)^{-1}{g}^2\preceq H^2$, from which we obtain
\begin{equation}
\begin{split}
\quad \int_D\left(2\mu  \epsilon({u_h}-{u_H}):\epsilon({u_h}-{u_H})+\lambda (\nabla\cdot({u_h}-{u_H}))^2\right) &\preceq
\big(\frac{1}{\Lambda_*}\big)^{m+1} \int_{D}\left( 2\mu  \epsilon({u_h}):\epsilon(u)+\lambda(\nabla\cdot{u_h})^2\right)\\
& \quad+\left((\Lambda_*)^m\left(\frac{1-(\Lambda_*)^{-m}}{\Lambda_*-1}\right)+1\right)H^2.
\end{split}
\end{equation}
Combining results above, we have
\begin{theorem}
Let $u\in V^h_{CG}$ be the fine-scale CG-FEM solution defined in (\ref{cg_fine_sol}) and $u_H$ be the  CG-GMsFEM solution defined in (\ref{cg_ms_sol}) without oversampling. If $\Lambda_*\ge1$ and $\int_D (\lambda+2\mu)^{-1}g^2\preceq1$, let $ n=-\frac{log(H)}{log\Lambda_*}$, then 
\begin{equation*}
\quad\int_D\left(2\mu  \epsilon({u_h}-{u_H}):\epsilon({u_h}-{u_H})+\lambda (\nabla\cdot({u_h}-{u_H}))^2\right)  \preceq\left(\frac{H}{\Lambda_*}\right)\left(\int_{D}\left(2\mu  \epsilon(u):\epsilon(u)+\lambda (\nabla\cdot{u_h})^2 \right)+1\right).
\end{equation*}
\end{theorem}

\subsection{Oversampling case}

In this subsection, we will analyze the convergence of CG-GMsFEM solution defined in (\ref{cg_ms_sol}) with oversampling.  We define $I^{\omega_n^{+}}{u_h}$ as an interpolation of 
${u_h}$ in $\omega_n^{+}$ using the first $L_n$ modes for the eigenvalue problem (\ref{cg_over_1}).
Let $\chi_n^+$ be a partition of unity subordinated to the coarse neighborhood $\omega_n^+$. We require $\chi_n^+$ to be zero on $\partial\omega_n^+$ and 
\begin{equation*}
|\nabla\chi_n|^2\preceq|\nabla\chi_n^+|^2.
\end{equation*}
Using the same argument as Lemma \ref{cacc_lemma}, it is easy to deduce
\begin{equation}
\begin{split}\label{cao}
&\quad\int_{\omega_n^+}\left(2\mu |\chi_n^+|^2\epsilon({u_h}-I^{\omega_n^+}{u_h}):\epsilon({u_h}-I^{\omega_n^+}{u_h})+\lambda |\chi_n^+|^2 (\nabla\cdot({u_h}-I^{\omega_n^+}{u_h}))^2 \right)\\
&\preceq |\int_{\omega_n^+} |\chi_n^+|^2{g}\cdot({u_h}-I^{\omega_n^+}{u_h})|+
\int_{\omega_n^+} (\lambda + 2\mu)|\nabla \chi_n^+|^2 ({u_h}-I^{\omega_n^+}{u_h})^2,
\end{split}
\end{equation}
where $g=f+\text{div}(\sigma(I^{\omega_n}{u_h}))$, $I^{\omega_n}{u_h} =I^{\omega_n^{+}}{u_h}$ in $\omega_n$.

Applying eigenvalue problem (\ref{cg_over_1}), we obtain
\begin{equation}\label{olocal}
\int_{\omega_n^+}(\lambda+2\mu) |\nabla\chi_n^+|^2 ({u_h}-I^{\omega_n^+}{u_h})^2\preceq \frac{1}{\xi_{L_n+1}^{\omega_n}}\int_{\omega_n}
\left(2\mu  \epsilon({u_h}-I^{\omega_n}{u_h}):\epsilon({u_h}-I^{\omega_n}{u_h}) +\lambda  (\nabla\cdot({u_h}-I^{\omega_n}{u_h}))^2\right).
\end{equation}
Using the definition of interpolation $I^{\omega_n^+}{u_h}$, we have
\begin{equation}
\begin{split}
&\quad\sum_{n=1}^{N_s}\int_{\omega_n}(\lambda+2\mu) |\nabla\chi_n|^2 ({u_h}-I^{\omega_n}{u_h})^2
\preceq\sum_{n=1}^{N_s}\int_{\omega_n^+}(\lambda+2\mu) |\nabla\chi_n^+|^2 ({u_h}-I^{\omega_n^+}{u_h})^2\\
&\preceq\sum_{n=1}^{N_s}  \frac{1}{\xi_{L_n+1}^{\omega_n}}\left(\int_{\omega_n}2\mu \epsilon({u_h}-I^{\omega_n}{u_h}):\epsilon({u_h}-I^{\omega_n}{u_h})+\int_{\omega_n}
\lambda  (\nabla\cdot({u_h}-I^{\omega_n}{u_h}))^2\right)   \\
&\preceq\sum_{n=1}^{N_s}  \frac{1}{\xi_{L_n+1}^{\omega_n}}\left(\int_{\omega_n^+} 2\mu  |\nabla\chi_n^+|^2\epsilon({u_h}-I^{\omega_n^+}{u_h}):\epsilon({u_h}-I^{\omega_n^+}{u_h})+\int_{\omega_n^+}
\lambda |\nabla\chi_n^+|^2(\nabla\cdot({u_h}-I^{\omega_n^+}{u_h}))^2\right) \\
&\preceq\sum_{n=1}^{N_s}\frac{1}{\xi_{L_n+1}^{\omega_n}}\int_{\omega_n^+} (\lambda + 2\mu)|\nabla \chi_n^+|^2 ({u_h}-I^{\omega_n^+}{u_h})^2+\sum_{n=1}^{N_s} \frac{1}{\xi_{L_n+1}^{\omega_n}}|\int_{\omega_n^+}|\chi_n^+|^2{g}\cdot({u_h}-I^{\omega_n^+}{u_h})|\\
&\preceq\frac{1}{\Lambda_*^+}\left(\sum_{n=1}^{N_s}\int_{\omega_n^+} (\lambda + 2\mu)|\nabla \chi_n^+|^2 ({u_h}-I^{\omega_n^+}{u_h})^2+\sum_{n=1}^{N_s} |\int_{\omega_n^+} |\chi_n^+|^2{g}\cdot({u_h}-I^{\omega_n^+}{u_h})|\right)\\
&\preceq\frac{1}{\Lambda_*^+}\sum_{n=1}^{N_s}\left(\frac{1}{\xi_{L_n+1}^{\omega_n}}\int_{\omega_n} \big(2\mu  \epsilon({u_h}-I^{\omega_n}{u_h}):\epsilon({u_h}-I^{\omega_n}{u_h})  +
\lambda  (\nabla\cdot({u_h}-I^{\omega_n}{u_h}))^2\big) + |\int_{\omega_n^+} |\chi_n^+|^2g\cdot({u_h}-I^{\omega_n^+}u)|\right),
\end{split}
\end{equation}
where ${\Lambda_*^+}=\text{min}_{\omega_n}\xi_{L_n+1}^{\omega_n}$.

 Applying the last inequality $m$ times with (\ref{olocal}),  we get
\begin{equation}
\begin{split}
&\quad\int_{\omega_n^+}(\lambda+2\mu) |\nabla\chi_n^+|^2 ({u_h}-I^{\omega_n^+}{u_h})^2\\
&\preceq \big(\frac{1}{\Lambda_*^+}\big)^m\left(\frac{1}{\xi_{L_n+1}^{\omega_n}}\int_{\omega_n} 2\mu \epsilon({u_h}-I^{\omega_n}{u_h}):\epsilon({u_h}-I^{\omega_n}{u_h}) + \frac{1}{\xi_{L_n+1}^{\omega_n}}\int_{\omega_n}\lambda  (\nabla\cdot({u_h}-I^{\omega_n}{u_h}))^2 \right)\\
&\quad+\sum_{l=1}^{m}\big(\frac{1}{\Lambda_*^+}\big)^l\sum_{n=1}^{N_s} |\int_{\omega_n^+}|\chi_n^+|^2 {g}\cdot({u_h}-I^{\omega_n^+}{u_h})|\\
&\preceq \big(\frac{1}{\Lambda_*^+}\big)^{m+1}\left(\int_{\omega_n} 2\mu \epsilon({u_h}-I^{\omega_n}{u_h}):\epsilon({u_h}-I^{\omega_n}{u_h})+\int_{\omega_n}\lambda (\nabla\cdot({u_h}-I^{\omega_n}{u_h}))^2 \right)\\
&\quad+(\Lambda_*^+)^m\left(\frac{1-(\Lambda_*^+)^{-m}}{\Lambda_*^+-1}\right)\sum_{n=1}^{N_s}\int_{\omega_n^+}((\lambda+2\mu)|\nabla\chi_n|^2)^{-1}{g}^2.\\
\end{split}
\end{equation}
Taking into account inequality (\ref{fem_estimate}), we have
\begin{equation}
\begin{split}
&\int_{D}\left(2\mu  \epsilon({u_h}-{u_H}):\epsilon({u_h}-{u_H})+\lambda (\nabla\cdot({u_h}-{u_H}))^2\right)\\
\preceq&\sum_{n=1}^{N_s} \int_{\omega_n}(\lambda +2\mu)|\nabla \chi_n|^2 ({u_h}-I^{\omega_n}{u_h})^2+\sum_{n=1}^{N_s}|\int_{\omega_n} \chi_n^2g\cdot({u_h}-I^{\omega_n}{u_h}) |\\
 \preceq& \big(\frac{1}{\Lambda_*^+}\big)^{m+1}\int_{D}   \left( 2\mu \epsilon (u):\epsilon(u)+\lambda (\nabla\cdot{u_h})^2\right)
+ \left((\Lambda_*^+)^m\left(\frac{1-(\Lambda_*^+)^{-m}}{\Lambda_*^+-1}\right)+1\right)R.
\end{split}
\end{equation}
where $R=\sum_{n=1}^{N_s}\int_{\omega_n}((\lambda+2\mu)|\nabla\chi_n^+|^2)^{-1}{g}^2$.

Therefore, similar with the no oversampling case, we have
\begin{theorem}
%$ {u_H}$ and $ {u_h}$ are CG-GMsFEM defined in (\ref{cg_ms_sol}) and fine-scale CG-FEM solution defined in (\ref{cg_fine_sol})
Let $u\in V^h_{CG}$ be the fine-scale CG-FEM solution defined in (\ref{cg_fine_sol}) and $u_H$ be the  CG-GMsFEM solution defined in (\ref{cg_ms_sol}) with oversampling. If $\Lambda_*^+\ge1$ and $\int_D (\lambda+2\mu)^{-1}g^2\preceq1$, let $ n=-\frac{log(H)}{log\Lambda_*^+}$, then 
\begin{equation*}
\quad\int_{D}\left(2\mu  \epsilon({u_h}-{u_H}):\epsilon({u_h}-{u_H})+\lambda (\nabla\cdot({u_h}-{u_H}))^2\right) \preceq\frac{H}{\Lambda_*^+}\left(\int_{D}  \left( 2\mu \epsilon (u):\epsilon(u)+\lambda (\nabla\cdot{u_h})^2\right)+1\right).
\end{equation*}
\end{theorem}
% For more the details about the estimation on $g$, we refer \cite{eglp13} which include the discussion for the elliptic equation.

\section{Error estimate for DG coupling} 
\label{sec:gmsfem_error_DG}

In this section, we will analyze the DG coupling of the GMsFEM (\ref{eq:ipdg}). 
For any $ {u}$, we define the DG-norm by
\begin{equation*}
\|  {u} \|_{\text{DG}}^2 = a_H( {u},  {u})
+ \sum_{E\in\mathcal{E}_{H}}\frac{\gamma}{h}\int_{E} \average{\lambda+2\mu}  \jump{{u}}^2 \; ds.
\end{equation*}
Let $K$ be a coarse grid block and let $ {n}_{\partial K}$ be the unit outward normal vector on $\partial K$.
We denote $V^h(\partial K)$
by the restriction of the conforming space $V^h$ on $\partial K$.
The normal flux $ {\sigma}( {u}) \,  {n}_{\partial K}$
is understood as an element in $V^h(\partial K)$ and is defined by
\begin{equation}
\int_{\partial K} ( {\sigma}( {u}) \,  {n}_{\partial K}) \cdot  {v} = 
\int_{K} \Big( 2\mu {\epsilon}( {u}): {\epsilon}(\widehat{ {v}}) 
+   \lambda \nabla\cdot  {u}  \nabla\cdot \widehat{ {v}}  \Big) \; d {x}, \quad  {v} \in V^h(\partial K),
\label{eq:flux}
\end{equation}
where $\widehat{ {v}}$ is the harmonic extension of $ {v}$ in $K$.
By the Cauchy-Schwarz inequality, 
\begin{equation*}
\int_{\partial K} ( {\sigma}( {u}) \,  {n}_{\partial K}) \cdot  {v} \leq 
a_H^K(u,u)^{\frac{1}{2}} \, a_H^K(\widehat{v},\widehat{v})^{\frac{1}{2}}.
\end{equation*}
By an inverse inequality and the fact that $\widehat{v}$ is the harmonic extension of $v$
\begin{equation*}
a_H^K(\widehat{v},\widehat{v}) \leq \kappa_K C^2_{\text{inv}} h^{-1} \int_{\partial K} |v|^2 \; dx,
\end{equation*}
where $\kappa_K = \max_K \{ \lambda+2\mu \}$ and $C_{\text{inv}} > 0$ is the constant from inverse inequality. 
Thus,
\begin{equation*}
\int_{\partial K} ( {\sigma}( {u}) \,  {n}_{\partial K}) \cdot  {v} \leq \kappa_K^{\frac{1}{2}} C_{\text{inv}} h^{-\frac{1}{2}} \| v\|_{L^2(\partial K)} \,
a_H^K(u,u)^{\frac{1}{2}}.
\end{equation*}
This shows that
\begin{equation*}
\int_{\partial K} | {\sigma}( {u}) \,  {n}_{\partial K} |^2
\leq  \kappa_K C^2_{\text{inv}} h^{-1} a_H^K(u,u).
\end{equation*}

Our first step in the convergence analysis 
is to establish the continuity and the coercivity 
of the bilinear form (\ref{eq:bilinear-ipdg})
with respect to the DG-norm. 

\begin{lemma}
Assume that the penalty parameter $\gamma$ is chosen so that $\gamma > 2 C_{\text{inv}}^2$. 
The bilinear form $a_{\text{DG}}$ defined in (\ref{eq:bilinear-ipdg})
is continuous and coercive, that is, 
\begin{eqnarray}
a_{\text{DG}}( {u},  {v})
&\leq&  \|  {u} \|_{\text{DG}} \, \|  {v} \|_{\text{DG}}, \\
a_{\text{DG}}( {u},  {u})
&\geq& a_0 \|  {u} \|_{\text{DG}}^2,
\end{eqnarray}
for all $ {u},  {v}$, where $a_0 = 1 - \sqrt{2} C_{\text{inv}} \gamma^{-\frac{1}{2}} >0$. 
\end{lemma}
{\it Proof}. By the definition of $a_{\text{DG}}$, we have
\begin{equation*}
a_{\text{DG}}( {u},  {v}) = a_H( {u},  {v})
- \sum_{E\in \mathcal{E}^H} \int_E \Big( \average{ {\sigma}( {u}) \,  {n}_E} \cdot \jump{ {v}} 
+ \average{ {\sigma}( {v}) \,  {n}_E} \cdot \jump{ {u}} \Big) \; ds
+ \sum_{E\in\mathcal{E}^H} \frac{\gamma}{h} \int_E \average{\lambda+2\mu} \jump{ {u}} \cdot \jump{ {v}} \; ds.
\end{equation*}
Notice that
\begin{equation*}
a_H( {u},  {v}) + \sum_{E\in\mathcal{E}^H} \frac{\gamma}{h} \int_E \average{\lambda+2\mu} \jump{ {u}} \cdot \jump{ {v}} \; ds
\leq  \| u\|_{\text{DG}} \, \| v\|_{\text{DG}}. 
\end{equation*}
For an interior coarse edge $E \in\mathcal{E}^H$, 
we let $K^{+}, K^{-}\in \mathcal{T}^H$ be the two coarse grid blocks having the edge $E$.
By the Cauchy-Schwarz ineqaulity, 
we have
\begin{equation}
\int_E  \average{ {\sigma}( {u}) \,  {n}_E} \cdot \jump{ {v}} \; ds
\leq \Big(  h\int_E  \average{ {\sigma}( {u}) \,  {n}_E}^2 \average{\lambda+2\mu}^{-1}  \; ds \Big)^{\frac{1}{2}}
\Big( \frac{1}{h} \int_E \average{\lambda+2\mu} \jump{ {v}}^2 \; ds \Big)^{\frac{1}{2}}.
\label{eq:cont1}
\end{equation}
Notice that
\begin{equation*}
\begin{split}
&\: h\int_E  \average{ {\sigma}( {u}) \,  {n}_E}^2 \average{\lambda+2\mu}^{-1}  \; ds \\
\leq &\: h\Big( \int_E ( {\sigma}( {u}^{+}) \,  {n}_E)^2 (\lambda^{+} + 2\mu^{+})^{-1} \; ds 
+ \int_E ( {\sigma}( {u}^{-}) \,  {n}_E)^2 (\lambda^{-} + 2\mu^{-})^{-1} \; ds 
\Big),
\end{split}
\end{equation*}
where $ {u}^{\pm} =  {u}|_{K^{\pm}}$, $\lambda^{\pm} = \lambda|_{K^{\pm}}$
and $\mu^{\pm} = \mu|_{K^{\pm}}$.
So, we have
\begin{equation*}
h\int_E  \average{ {\sigma}( {u}) \,  {n}_E}^2 \average{\lambda+2\mu}^{-1}  \; ds
\leq C_{\text{inv}}^2 \Big( a_H^{K^{+}}(u^{+},u^{+}) + a_H^{K^{-}}(u^{-},u^{-}) \Big). 
\end{equation*}
Thus (\ref{eq:cont1}) becomes
\begin{equation}
\int_E  \average{ {\sigma}( {u}) \,  {n}_E} \cdot \jump{ {v}} \; ds \leq 
C_{\text{inv}} \Big( a_H^{K^{+}}(u^{+},u^{+}) + a_H^{K^{-}}(u^{-},u^{-}) \Big)^{\frac{1}{2}}
\Big( \frac{1}{h} \int_E \average{\lambda+2\mu} \jump{ {v}}^2 \; ds \Big)^{\frac{1}{2}}.
\label{eq:cont2}
\end{equation}
When $E$ is a boundary edge, we have
\begin{equation}
\int_E  \average{ {\sigma}( {u}) \,  {n}_E} \cdot \jump{ {v}} \; ds \leq 
C_{\text{inv}}  a_H^{K}(u,u)^{\frac{1}{2}}
\Big( \frac{1}{h} \int_E \average{\lambda+2\mu} \jump{ {v}}^2 \; ds \Big)^{\frac{1}{2}},
\label{eq:cont3}
\end{equation}
where $K$ denotes the coarse grid block having the edge $E$.
Summing (\ref{eq:cont2}) and (\ref{eq:cont3}) for all edges $E\in\mathcal{E}^H$, 
we have
\begin{equation*}
\sum_{E\in \mathcal{E}^H} \int_E  \average{ {\sigma}( {u}) \,  {n}_E} \cdot \jump{ {v}} \; ds
\leq \sqrt{2} C_{\text{inv}} a_H(u,u)^{\frac{1}{2}} \Big( \sum_{E\in \mathcal{E}^H} \frac{1}{h} \int_E \average{\lambda+2\mu} \jump{ {v}}^2 \; ds \Big)^{\frac{1}{2}}.
\end{equation*}
Similarly, we have
\begin{equation*}
\sum_{E\in \mathcal{E}^H} \int_E  \average{ {\sigma}( {v}) \,  {n}_E} \cdot \jump{ {u}} \; ds
\leq \sqrt{2} C_{\text{inv}} a_H(v,v)^{\frac{1}{2}} \Big( \sum_{E\in \mathcal{E}^H} \frac{1}{h} \int_E \average{\lambda+2\mu} \jump{ {u}}^2 \; ds \Big)^{\frac{1}{2}}.
\end{equation*}
Hence
\begin{equation}
\sum_{E\in \mathcal{E}^H} \int_E \Big( \average{ {\sigma}( {u}) \,  {n}_E} \cdot \jump{ {v}} 
+ \average{ {\sigma}( {v}) \,  {n}_E} \cdot \jump{ {u}} \Big) \; ds
\leq \sqrt{2} C_{\text{inv}} \gamma^{-\frac{1}{2}} \| u\|_{\text{DG}} \, \|v\|_{\text{DG}}.
\label{eq:cont4}
\end{equation}
This proves the continuity. 

For coercivity, we have
\begin{equation*}
a_{\text{DG}}( {u},  {u}) = \| u\|_{\text{DG}}^2
- \sum_{E\in \mathcal{E}^H} \int_E \Big( \average{ {\sigma}( {u}) \,  {n}_E} \cdot \jump{ {u}} 
+ \average{ {\sigma}( {u}) \,  {n}_E} \cdot \jump{ {u}} \Big) \; ds.
\end{equation*}
By (\ref{eq:cont4}), we have
\begin{equation*}
a_{\text{DG}}( {u},  {u}) \geq (1 - \sqrt{2} C_{\text{inv}} \gamma^{-\frac{1}{2}} ) \| u\|_{\text{DG}}^2,
\end{equation*}
which gives the desired result. 
\begin{flushright}
$\square$
\end{flushright}

We will now prove the convergence of the method (\ref{eq:ipdg}). 
Let $u_h \in V^h_{\text{DG}}$ be the fine grid solution which satisfies
\begin{equation}
a_{\text{DG}}(u_h, v) = (f,v), \quad\forall v\in V^h_{\text{DG}}.
\label{eq:ipdgfine}
\end{equation}
It is well-known that $u_h$ converges to the exact solution $u$ in the DG-norm as the fine mesh size $h\rightarrow 0$. 
Next, we define a projection $u_S\in V^{\text{snap}}$ of $u_h$ in the snapshot space 
by the following construction. 
For each coarse grid block $K$, the restriction of $u_S$ on $K$
is defined as the harmonic extension of $u_h$, that is,
\begin{equation}
\begin{split}
- \nabla \cdot  {\sigma}( {u}_{S}) &=  {0}, \quad\text{ in } \; K, \\
 {u}_{S} &= u_h, \quad\text{ on } \; \partial K.
\end{split}
\label{eq:proj}
\end{equation}
Now, we prove the following estimate for the projection $u_S$. 

\begin{lemma}
Let $u_h\in V^h_{\text{DG}}$ be the fine grid solution defined in (\ref{eq:ipdgfine})
and $u_S\in V^{\text{snap}}$ be the projection of $u_h$ defined in (\ref{eq:proj}).
Then we have
\begin{equation*}
\| u_h - u_S \|_{\text{DG}} \leq CH \Big( \max_{K\in\mathcal{T}^H} \eta_K \Big) \| f \|_{L^2(\Omega)},
\end{equation*}
where $\eta_K = \min_K \{\lambda+2\mu\}$. 
\end{lemma}
{\it Proof}. Let $K$ be a given coarse grid block.
Since $u_S=u_h$ on $\partial K$,
the jump terms in the DG-norm vanish. 
Thus, the DG-norm can be written as
\begin{equation*}
\| u_h - u_S \|_{\text{DG}}^2 = \sum_{K\in\mathcal{T}^H} a_H^K(u_h-u_S,u_h-u_S).
\end{equation*}
Since $u_S$ satisfies (\ref{eq:proj}) and $u_h-u_S=0$ on $\partial K$, we have
\begin{equation*}
a_H^K(u_S, u_h-u_S) = 0. 
\end{equation*}
So,
\begin{equation*}
\| u_h - u_S \|_{\text{DG}}^2 = \sum_{K\in\mathcal{T}^H} a_H^K(u_h,u_h-u_S)
= a_{\text{DG}}(u_h, u_h-u_S)
= (f, u_h-u_S).
\end{equation*}
By the Poincare inequality, we have
\begin{equation*}
\| u_h - u_S \|_{L^2(K)} \leq C H^2 \eta^2_K  a_H^K(u_h-u_S,u_h-u_S),
\end{equation*}
where $\eta_K = \min_K \{ \lambda+2\mu \}$.
Hence, we have
\begin{equation*}
\| u_h - u_S \|_{\text{DG}} \leq CH \Big( \max_{K\in\mathcal{T}^H} \eta_K \Big) \| f \|_{L^2(\Omega)}. 
\end{equation*} 

\begin{flushright}
$\square$
\end{flushright}

In the following theorem, we will state and prove the convergence 
of the GMsFEM (\ref{eq:ipdg}).

\begin{theorem}
Let $u_h\in V^h_{\text{DG}}$ be the fine grid solution defined in (\ref{eq:ipdgfine})
and $u_H$ be the GMsFEM solution defined in (\ref{eq:ipdg}). Then we have
\begin{equation*}
\| u_h - u_H \|_{\text{DG}}^2 
\leq C\Big( \sum_{i=1}^{N_E} \frac{H}{\left\langle\lambda+2\mu\right\rangle\xi_{L_i+1}} ( 1 + \frac{\gamma H}{h \xi_{L_i+1}} ) \int_{\partial K_i} ( \sigma(u_S)\cdot n_{\partial K})^2 \; ds + 
H^2 \Big( \max_{K\in\mathcal{T}^H} \eta_K^2 \Big) \| f \|^2_{L^2(\Omega)} \Big),
\end{equation*}
where $u_S$ is defined in (\ref{eq:proj}). 
\end{theorem}
{\it Proof}.
First, we will define a projection $\widehat{u}_S \in V^{\text{off}}$
of $u_S$ in the offline space. 
Notice that, on each $K_i$, $u_S$ can be represented by
\begin{equation*}
u_S = \sum_{l=1}^{M_i} c_l \psi_l^{i,\text{off}},
\end{equation*}
where $M_i = M^{i,\text{snap}}$ and 
we assume that the functions $\psi_l^{i,\text{off}}$
are normalized so that
\begin{equation*}
 \int_{\partial K_i} \left\langle\lambda+2\mu\right\rangle (\psi_l^{i,\text{off}})^2 \; ds = 1.
\end{equation*}
Then the function $\widehat{u}_S$ is defined by
\begin{equation*}
\widehat{u}_S = \sum_{l=1}^{L_i} c_l \psi_l^{i,\text{off}}.
\end{equation*}
We will find an estimate of $\| u_S - \widehat{u}_S\|_{\text{DG}}$. 
Let $K$ be a given coarse grid block. 
Recall that the spectral problem is
\begin{equation*}
\int_{K} 2\mu {\epsilon}({u}): {\epsilon}({v})dx 
+ \int_K  \lambda \nabla\cdot  {u}  \nabla\cdot {v}
= \frac{\xi}{H} \int_{\partial K} \left\langle\lambda+2\mu\right\rangle {u} {v} \; ds.
\end{equation*}
By the definition of the flux (\ref{eq:flux}), the spectral problem can be represented as
\begin{equation*}
\int_{\partial K} (\sigma(u)\cdot n_{\partial K}) v \; ds= \frac{\xi}{H} \int_{\partial K}\left\langle\lambda+2\mu\right\rangle {u} {v} \; ds.
\end{equation*}
By the definition of the DG-norm, the error $\| u_S - \widehat{u}_S\|_{\text{DG}}$ can be computed as
\begin{equation*}
\| \widehat{{u}}_S - {u}_S \|_{\text{DG}}^2
\leq \sum_{K} \Big( \int_K 2\mu  { \epsilon}(\widehat{{u}}_S - {u}_S )^2 
+ \int_K \lambda (\nabla \cdot (\widehat{{u}}_S - {u}_S) )^2 
+  \frac{\gamma}{h} \int_{\partial K} \average{\lambda+2\mu} (\widehat{{u}}_S - {u}_S))^2 \Big).
\end{equation*}
Note that
\begin{equation*}
\int_{K_i} 2\mu  { \epsilon}(\widehat{{u}}_S - {u}_S )^2 
+ \int_{K_i} \lambda (\nabla \cdot (\widehat{{u}}_S - {u}_S) )^2
\leq\frac{1}{h}\int_{\partial{K_i}}\left\langle\lambda+2\mu\right\rangle (\widehat{{u}}_S - {u}_S)^2
= \sum_{l={L_i}+1}^{M_i}  \frac{\xi_l}{H} c_l^2
\leq \frac{H}{\xi_{L_i+1}}  \sum_{l=L_i+1}^{M_i} (\frac{\xi_l}{H})^2 c_l^2.
\end{equation*}
Also,
\begin{equation*}
\frac{1}{h} \int_{\partial K_i} \average{\lambda+2\mu} (\widehat{{u}}_S - {u}_S)^2 
= \frac{1}{h} \sum_{l=L_i+1}^{M_i} c_l^2
\leq \frac{H^2}{h \xi_{L_i+1}^2} \sum_{l=L_i+1}^{M_i} (\frac{\xi_l}{H})^2 c_l^2.
\end{equation*}
Moreover, 
\begin{equation*}
\sum_{l=L_i+1}^{M_i} (\frac{\xi_l}{H})^2 c_l^2 \leq \sum_{l=1}^{M_i} (\frac{\xi_l}{H})^2 c_l^2 \leq \frac{1}{\left\langle\lambda+2\mu\right\rangle}\int_{\partial K_i} ( \sigma(u_S)\cdot n_{\partial K})^2 \; ds.
\end{equation*}
Consequently, we obtain the following bound
\begin{equation*}
\| u_S - \widehat{{u}}_S  \|^2_{\text{DG}}
\leq \sum_{i=1}^{N_E} \frac{H}{\left\langle\lambda+2\mu\right\rangle\xi_{L_i+1}} ( 1 + \frac{\gamma H}{h \xi_{L_i+1}} ) \int_{\partial K_i} ( \sigma(u_S)\cdot n_{\partial K})^2 \; ds.
\end{equation*}

Next, we will prove the required error bound.
By coercivity,
\begin{equation*}
\begin{split}
a_0 \| \widehat{u}_S - u_H \|_{\text{DG}}^2 
&= a_{\text{DG}}(\widehat{u}_S-u_H, \widehat{u}_S-u_H) \\
&= a_{\text{DG}}(\widehat{u}_S-u_H, \widehat{u}_S-u_S)   + a_{\text{DG}}(\widehat{u}_S-u_H, u_S-u_h) + a_{\text{DG}}(\widehat{u}_S-u_H, u_h-u_H).
\end{split}
\end{equation*}
Note that $a_{\text{DG}}(\widehat{u}_S-u_H, u_h-u_H) = 0$ since $\widehat{u}-u_H \in V^{\text{off}}$. 
Using the above results,
\begin{equation*}
\| \widehat{u}_S - u_H \|_{\text{DG}}^2 
\leq C\Big( \sum_{i=1}^{N_E} \frac{H}{\left\langle\lambda+2\mu\right\rangle\xi_{L_i+1}} ( 1 + \frac{\gamma H}{h \xi_{L_i+1}} ) \int_{\partial K_i} ( \sigma(u_S)\cdot n_{\partial K})^2 \; ds + 
H^2 \Big( \max_{K\in\mathcal{T}^H} \eta_K^2 \Big) \| f \|^2_{L^2(\Omega)} \Big).
\end{equation*}
Finally, the desired bound is obtained by the triangle inequality
\begin{equation*}
\| u_h - u_H\|_{\text{DG}} \leq \| u_h - u_S\|_{\text{DG}} + \| u_S - \widehat{u}_S\|_{\text{DG}} + \| \widehat{u}_S - u_H\|_{\text{DG}}.
\end{equation*}

\begin{flushright}
$\square$
\end{flushright}

\section{Conclusions}

In this paper, we design a multiscale model reduction method using GMsFEM
for elasticity equations in heterogeneous media. 
%Our applications are motivated 
%elastic wave propagation in subsurface where subsurface properties 
%contain multiple scales and high contrast. 
We design a snapshot space and
an offline space based on the analysis. We present two approaches that couple
multiscale basis functions of the offline space. These are continuous Galerkin
and discontinuous Galerkin methods. Both approaches are analyzed. 
We present oversampling studies where larger domains are used for
calculating the snapshot space.
Numerical
results are presented.

\bibliographystyle{plain}   
\bibliography{references1}

\end{document}